\documentclass[a4paper]{amsart}
\usepackage{amssymb} 
\usepackage[T1]{fontenc}
\usepackage[latin1]{inputenc}
\usepackage{amsfonts}
\usepackage{amsxtra}
\usepackage{ae}
\usepackage[all]{xy}
\usepackage{enumerate}

\include{diagram}

\newcommand*{\ket}{\rangle}
\newcommand*{\bra}{\langle}

\newcommand*{\Smooth}{\mathfrak{Smooth}}

\newcommand*{\C}{\mathcal{C}}
\newcommand*{\D}{\mathcal{D}}
\newcommand*{\E}{\mathcal{E}}

\newcommand*{\ayd}{\mathsf{AYD}}
\newcommand*{\AYD}{\mathsf{A}}
\newcommand*{\SHom}{\mathfrak{Hom}}

\newcommand*{\an}{\textsf{an}}
\newcommand*{\Or}{\textsf{Or}}

\newcommand*{\LH}{\mathbb{L}}

\newcommand*{\red}{\mathsf{red}}

\DeclareMathOperator{\ind}{ind}
\DeclareMathOperator{\res}{res}
\DeclareMathOperator{\Aut}{Aut}

\DeclareMathOperator{\Hom}{Hom}

\DeclareMathOperator{\id}{id}

\DeclareMathOperator{\ch}{ch}

\DeclareMathOperator{\hyp}{hyp}
\DeclareMathOperator{\elli}{ell}

\DeclareMathOperator{\tr}{tr}

\numberwithin{equation}{section}
\theoremstyle{change}
\newtheorem{theorem}{Theorem}[section]
\newtheorem{prop}[theorem]{Proposition}
\newtheorem{lemma}[theorem]{Lemma}

\newtheorem{definition}[theorem]{Definition}

\begin{document}

\title[Chern character for totally disconnected groups]{Chern character for totally disconnected groups}
\author{Christian Voigt}
\address{Institut for Mathematical Sciences\\
         University of Copenhagen\\
         Universitetsparken 5 \\
         2100 Copenhagen\\
         Denmark
}
\email{cvoigt@math.ku.dk}

\subjclass[2000]{46L80, 19D55, 55N91, 19L47}

\maketitle

\begin{abstract}
In this paper we construct a bivariant Chern character for the equivariant $ KK $-theory 
of a totally disconnected group with values in bivariant equivariant cohomology in the sense of 
Baum and Schneider. We prove in particular that the complexified left hand side of the Baum-Connes 
conjecture for a totally disconnected group is isomorphic to cosheaf homology. 
Moreover, it is shown that our transformation extends the Chern character defined by Baum and 
Schneider for profinite groups. 
\end{abstract}

\section{Introduction}

Let $ G $ be a second countable locally compact group. 
The Baum-Connes conjecture \cite{BCH} asserts that the $ K $-theory of the reduced group 
$ C^* $-algebra $ C^*_\red(G) $ of $ G $ is isomorphic to the 
equivariant $ K $-homology of the universal proper $ G $-space $ \underline{E}G $. More precisely, 
the assembly map
$$ 
\mu: K^G_*(\underline{E}G) \rightarrow K_*(C^*_\red(G))
$$ 
is conjectured to be an isomorphism. 
In this paper we construct an equivariant Chern character for the left hand side of the assembly map in the case that the group $ G $ is 
totally disconnected. 
We show in particular that the complexified equivariant $ K $-homology of $ \underline{E}G $ for a totally disconnected group $ G $ 
is given by cosheaf homology \cite{BCH}. \\
The corresponding result for discrete groups has been obtained independently by Baum and Connes \cite{BC2} and L\"uck \cite{Lueck}. 
The approach of L\"uck actually yields Chern characters for arbitrary equivariant homology theories 
with source given by associated Bredon homologies. For equivariant $ K $-homology and group actions on simplicial 
complexes, the corresponding Bredon homology with complex coefficients is naturally isomorphic to cosheaf homology. 
However, it seems to be unclear wether the case of totally disconnected groups 
can still be handled in the framework of L\"uck. We remark that the result for totally disconnected groups stated above has also been 
obtained by Baum, Block and Higson in an unpublished paper \cite{BBH}. \\
In fact we prove a more general statement. The target of our Chern character is bivariant 
equivariant cohomology in the sense of Baum and Schneider \cite{BS}. One virtue of this
cohomology theory is that it unifies and generalizes previous constructions in the literature. In 
particular, the theory of Baum and Schneider contains the cosheaf homology groups mentioned 
above as a special case. \\
The main result of this paper is the construction of a bivariant Chern character 
\begin{equation*}
ch_*^G: KK^G_*(C_0(X),C_0(Y)) \rightarrow \bigoplus_{j \in \mathbb{Z}} H^{* + 2j}_G(X,Y)
\end{equation*}
from equivariant $ KK $-theory with values in the theory of Baum and Schneider. 
Here $ G $ is a totally disconnected group, $ X $ and $ Y $ are finite dimensional locally finite $ G $-simplicial complexes in the sense of \cite{Voigtbs} 
and $ X $ is assumed to be proper and $ G $-finite. This Chern character becomes an isomorphism after tensoring the left hand side
with the complex numbers. We also discuss how the character can be extended to $ G $-$ CW $-complexes. The 
restriction to simplicial complexes is necessary in order to apply the machinery of equivariant cyclic homology \cite{Voigtepch}.  
In fact, equivariant cylic homology is the main ingredient in our construction. \\
We remark that Baum and Schneider conjectured the existence of such a bivariant Chern character in \cite{BS} and 
proved this conjecture in the case of profinite groups. However, their construction does not extend to noncompact groups. 
Thus, in some sense, the present paper completes the work of Baum and Schneider. \\
An important ingredient in the construction of L\"uck is the induction structure for an equivariant homology 
theory \cite{Lueck}. As already mentioned above, the approach presented here is based on equivariant cyclic homology. 
In turn, for the construction of our Chern character we need some information on the 
compatibility of equivariant cyclic homology with induction. From a technical point of view the corresponding 
considerations constitute the main part of this paper. Apart from this, the constructions and results from \cite{Voigtelch} 
play a central role. As a matter of fact, our Chern character is obtained by combining the 
equivariant Chern-Connes character from \cite{Voigtelch} with the description of bivariant equivariant 
cohomology in terms of cyclic homology given in \cite{Voigtbs}. \\
Let us explain how the paper is organized. In section \ref{secsmoothcov} we recall some facts 
about smooth representations of totally disconnected groups and anti-Yetter-Drinfeld modules. 
Smooth representations and anti-Yetter-Drinfeld modules are a basic ingredient in the definition of equivariant cyclic 
homology and the theory of Baum and Schneider. Section \ref{secech} contains a brief 
survey of the different variants of equivariant cyclic homology involved in the construction of the Chern character. 
The basic theory of restriction and induction of $ G $-modules, $ G $-algebras and anti-Yetter-Drinfeld modules 
is discussed in section \ref{secindmodules}. In section \ref{secindcyclic} this is used to study induction in 
equivariant cyclic homology. 
Section \ref{secchern} contains the construction of the equivariant Chern character for $ G $-simplicial complexes. 
This yields the main result stated above. In section \ref{secbs} we show 
that our Chern character restricts to the natural transformation defined by Baum and Schneider in the case of profinite groups. 
Finally, in section \ref{seccw} we explain how to extend the homological Chern character arising from our constructions 
to proper $ G $-$ CW $-complexes. In the setting of discrete groups this shows in particular that our transformation may
be defined on the same class of spaces as the Chern character constructed by L\"uck. However, we have not studied the relation 
between these characters. \\
I would like to thank W. L\"uck for some helpful comments. 

\section{Smooth representations and anti-Yetter-Drinfeld modules} \label{secsmoothcov}

Throughout this paper let $ G $ be a second countable locally compact and totally disconnected group. 
As in \cite{Voigtbs} we call 
an element $ t \in G $ elliptic if it is contained in a compact subgroup. The set of all elliptic elements 
of $ G $ is denoted by $ G_{\elli} $. In contrast 
we shall say that an element $ t \in G $ is hyperbolic if it is not elliptic. Let $ G_{\hyp} $ be the set 
of all hyperbolic elements of $ G $. Hence, according to these definitions, we obtain a disjoint union decomposition 
$$ 
G = G_{\elli} \cup G_{\hyp} 
$$ 
of the space $ G $. It follows from the structure theory developped in \cite{Willis} that $ G_{\elli} $ is a closed subset of $ G $ and 
that $ G_{\elli} $ is the union of all compact open subgroups of $ G $. Consequently $ G_{\hyp} $ is again an open and closed 
subset of $ G $. \\
In this section we review some facts about smooth representations of totally disconnected groups 
on bornological vector spaces and anti-Yetter-Drinfeld modules. For more information on bornological vector spaces and smooth 
representations we refer to \cite{H-L2}, \cite{HLMosc}, \cite{Meyerthesis}, \cite{Meyersmoothrep}. 
Concerning anti-Yetter-Drinfeld modules further details can be found in \cite{Voigtepch}, \cite{Voigtbs} where these objects were called covariant modules instead. 
All bornological vector spaces in this paper are assumed to be convex and separated. The completion of a bornological vector space $ V $ is 
denoted by $ V^c $. \\
A representation of $ G $ on a bornological vector space $ V $ is a group homomorphism 
$ \pi: G \rightarrow \Aut(V) $ where $ \Aut(V) $ denotes the group of bounded linear automorphisms of $ V $. 
A bounded linear map $ f: V \rightarrow W $ between representations of  $ G $ 
is called equivariant if $ f(s \cdot v) = s \cdot f(v) $ for all $ s \in G, v \in V $. 
To a representation $ \pi: G \rightarrow \Aut(V) $ one associates the linear map $ [\pi]: V \rightarrow F(G,V) $ defined by 
$ [\pi](v)(t) = \pi(t)(v) $ where $ F(G,V) $ denotes the space of all functions from 
$ G $ to $ V $. We denote by $ \E(G,V) $ the space of all locally constant functions on $ G $ with 
values in $ V $. 
\begin{definition} Let $ V $ be a separated bornological vector space. 
A representation $ \pi $ of $ G $ on $ V $ is called smooth if $ [\pi] $ defines a bounded linear map from 
$ V $ into $ \E(G,V) $. A separated (complete) $ G $-module is a smooth representation of $ V $ 
on a separated (complete) bornological vector space. 
\end{definition}
If $ V $ is a separated $ G $-module the stabilizers of small subsets of $ V $ are open subgroups of $ G $. 
Conversely, if $ \pi $ is a representation on a bornological vector space $ V $ such that the stablizers of all small subsets are open it follows 
that $ \pi $ is a $ G $-module. In particular, if $ V $ carries the fine bornology one recovers the usual notion of a smooth 
representation on a complex vector space. \\
The Hecke algebra of the group $ G $ is the space $ \D(G) $ of locally constant functions with 
compact support on $ G $ equipped with the convolution product with respect to a fixed 
Haar measure. Every separated $ G $-module $ V $ becomes a module over $ \D(G) $ using integration. Moreover, the $ \D(G) $-module $ V $ 
obtained in this way is essential in the sense that the canonical map $ \D(G) \otimes_{\D(G)} V \rightarrow V $ is an isomorphism. 
Conversely, every essential $ \D(G) $-module is the integrated form of a smooth representation. \\
For an arbitrary representation of $ G $ on a separated bornological vector space $ V $ there exists a
smoothing $ \Smooth(V) $ which is a smooth representation. 
This construction yields a functor $ \Smooth $ that is right adjoint to the natural forgetful functor from the 
category of $ G $-modules to the category of arbitrary representations on bornological vector spaces. \\
Let us now recall the concept of an anti-Yetter-Drinfeld module. We write $ \mathcal{O}_G $ for the vector space 
$ \D(G) $ equipped with pointwise multiplication and the action of $ G $ by conjugation. 
\begin{definition}
A separated (complete) anti-Yetter-Drinfeld module $ M $ is a separated (complete) $ G $-module which is at the same time 
an essential $ \mathcal{O}_G $-module. The $ G $-module structure and the $ \mathcal{O}_G $-module structure 
are required to be compatible in the sense that 
$$
s \cdot (f \cdot m) = (s \cdot f) \cdot (f \cdot m)
$$
for all $ s \in G, f \in \mathcal{O}_G $ and $ m \in M $. 
\end{definition}
A homomorphism $ \xi: M \rightarrow N $ of anti-Yetter-Drinfeld modules 
is a bounded linear map which is both $ G $-equivariant and $ \mathcal{O}_G $-linear. 
In the sequel we will use the abbreviations $ \ayd $-module and $ \ayd $-map for 
anti-Yetter-Drinfeld modules and their homomorphisms. We write $ \SHom_G(M, N) $ for the space of 
$ \ayd $-maps between $ \ayd $-modules $ M $ and $ N $. \\
The category of $ \ayd $-modules is isomorphic to the category of essential $ \AYD(G) $-modules where 
$ \AYD(G) = \mathcal{O}_G \rtimes G $ is the smooth crossed product of $ \mathcal{O}_G $ with respect to the adjoint action. 
In particular, for any bornological vector space $ V $ the space $ \AYD(G) \otimes V $ becomes an $ \ayd $-module. 
The modules arising in this way are projective, and every projective $ \ayd $-module is a direct summand in 
an $ \ayd $-module of this form. 
An important property of $ \ayd $-modules is 
the fact that they are equipped with a natural automorphism $ T $. In this way the category of $ \ayd $-modules 
becomes a para-additive category in the sense of \cite{Voigtepch}. For equivariant differential forms we will describe the automorphism $ T $ 
explicitly in section \ref{secech} below. \\
Let us discuss some features of $ \ayd $-modules which are specific to the setting of totally disconnected groups. 
For a totally disconnected group there is a natural decomposition of an $ \ayd $-module into an elliptic and 
a hyperbolic part. Since the group $ G $ is the disjoint union of its elliptic and hyperbolic elements one has   
two natural multipliers $ P_{\elli} $ and $ P_{\hyp} $ of the algebra $ \mathcal{O}_G $. By definition, the multiplier $ P_{\elli} $ is 
the characteristic function of the set $ G_{\elli} $ whereas $ P_{\hyp} = 1 - P_{\elli} $ is the characteristic 
function of $ G_{\hyp} $. For every $ \ayd $-module $ M $ this yields a natural direct sum decomposition of $ \ayd $-modules 
$$ 
M = M_{\elli} \oplus M_{\hyp} 
$$ 
where $ M_{\elli} = P_{\elli} \cdot M $ and $ M_{\hyp} = P_{\hyp} \cdot M $. Usually, the behaviour of the operator $ T $ 
is quite different on the two summands in this decomposition. For instance, for every element $ m \in M_{\elli} $ there 
exists a natural number $ n $ such that $ T^n(m) = m $. Using this fact one constructs 
an $ \ayd $-map $ E: M_{\elli} \rightarrow M_{\elli} $ which yields a projection onto the $ T $-invariant 
elements in $ M_{\elli} $. In contrast, if $ M $ is a projective $ \ayd $-module then the operator $ \id - T $ is injective 
on the hyperbolic part $ M_{\hyp} $ of $ M $. Moreover $ \id - T:M_{\hyp} \rightarrow (\id - T)M_{\hyp} $ 
is an isomorphism in this case. \\
If $ M $ is an $ \ayd $-module we write $ M \otimes _T \mathbb{C} $ for the space of $ T $-coinvariants, that is, the quotient of $ M $ 
by the closed subspace generated by all elements $ T(m) - m $ for $ m \in M $. Note that the action of $ T $ on $ M $ can be viewed as a module 
structure of the algebra $ R = \mathbb{C}[z,z^{-1}] $ of Laurent polynomials. Then $ M \otimes_T \mathbb{C} = M \otimes_R \mathbb{C}  $ where 
$ \mathbb{C} $ is viewed as a module over $ R $ using the character given by evaluation at $ 1 $. Remark that $ M \otimes_T \mathbb{C} $ 
is again an $ \ayd $-module in a natural way. Every $ \ayd $-map $ M \rightarrow N $ induces an $ \ayd $-map 
$ M \otimes_T \mathbb{C}\rightarrow N \otimes_T \mathbb{C} $ on $ T $-coinvariants. \\
Finally, for local cyclic homology the concept of a primitive $ \ayd $-module is important \cite{Voigtelch}. An $ \ayd $-module 
$ P $ is called primitive if there is a small disk $ S \subset P $ such that the natural map 
$ \AYD(G) \otimes \bra S \ket \rightarrow P $ is a quotient map. Every $ \ayd $-module $ M $ can be 
written in a canonical way as direct limit of the primitive modules $ \AYD(G)\bra S \ket $ generated by small disks $ S $ 
in $ M $. The relation between inductive systems of primitive modules and general $ \ayd $-modules is analogous to the one between 
inductive systems of normed spaces and separated bornological vector spaces. 

\section{Equivariant periodic, analytic and local cyclic homology}\label{secech}

In this section we review the definition of equivariant periodic, analytic and local cyclic homology. 
These theories are studied in detail in \cite{Voigtepch}, \cite{Voigtelch}. \\
A separated $ G $-algebra is a separated bornological algebra $ A $ which is also a $ G $-module 
such that the multiplication map $ A \otimes A \rightarrow A $ is equivariant. The unitarization of 
a $ G $-algebra $ A $ is denoted by $ A^+ $ and again a $ G $-algebra in a natural way where the group acts trivially 
on scalar multiples of the unit element $ 1 \in A^+ $. \\
Let $ A $ be a separated $ G $-algebra. The equivariant $ n $-forms of $ A $ are defined by 
$ \Omega^n_G(A) = \mathcal{O}_G \otimes \Omega^n(A) $ where $ \Omega^n(A) = A^+ \otimes A^{\otimes n} $ is 
the space of noncommutative differential forms over $ A $. The group $ G $ acts diagonally on $ \Omega^n_G(A) $ and there is
an obvious $ \mathcal{O}_G $-module structure given by multiplication. In this way $ \Omega^n_G(A) $ becomes a separated $ \ayd $-module. 
We write $ \Omega_G(A) $ for the direct sum of the spaces $ \Omega^n_G(A) $. \\
On equivariant differential forms there are two important boundary operators. The equivariant Hochschild 
boundary $ b: \Omega^n_G(A) \rightarrow \Omega^{n - 1}_G(A) $ is defined by 
\begin{align*}
b(f(s) \otimes &x_0 dx_1 \cdots dx_n) = f(t) \otimes x_0 x_1 dx_2 \cdots dx_n \\
&+ \sum_{j = 1}^{n - 1} (-1)^j f(s) \otimes x_0 dx_1 \cdots d(x_j x_{j + 1}) \cdots dx_n \\ 
&+ (-1)^n f(s) \otimes (s^{-1} \cdot x_n)x_0 dx_1 \cdots dx_{n - 1}.
\end{align*}
Moreover we have the equivariant Connes operator $ B: \Omega^n_G(A) \rightarrow \Omega^{n + 1}_G(A) $ which is given by 
\begin{equation*}
B(f(s) \otimes x_0dx_1 \cdots dx_n) = \sum_{i = 0}^n (-1)^{ni}
f(s) \otimes s^{-1} \cdot(dx_{n + 1 - i} \cdots dx_n)dx_0 \cdots
dx_{n - i}.
\end{equation*}
Both $ b $ and $ B $ are $ \ayd $-maps and the natural isomorphism $ T $ for 
$ \ayd $-modules has the form
\begin{equation*}
T(f(s) \otimes \omega) = f(s) \otimes s^{-1} \cdot \omega 
\end{equation*}
on $ \Omega^n_G(A) $. 
The space of equivariant differential forms together with the operators $ b $ and $ B $ forms a paramixed complex \cite{Voigtepch}
which means that the relations $ b^2 = 0 $, $ B^2 = 0 $ and
$ [b,B] = bB + Bb = \id - T $ hold. \\
Equivariant differential forms are used to define the equivariant $ X $-complex of a separated $ G $-algebra. 
\begin{definition} Let $ A $ be a separated $ G $-algebra. The equivariant $ X $-complex 
$ X_G(A) $ of $ A $ is 
\begin{equation*}
    \xymatrix{
      {X_G(A) \colon \ }
        {\Omega^0_G(A)\;} \ar@<1ex>@{->}[r]^{B} &
          {\;\Omega^1_G(A)/ b(\Omega^2_G(A)).} 
            \ar@<1ex>@{->}[l]^-{b} 
               } 
\end{equation*}
\end{definition}
If $ \partial $ denotes the boundary operator in $ X_G(A) $ then $ \partial^2 = \id - T $ is 
in general not zero. Hence the equivariant $ X $-complex is not a chain complex. 
Instead, it is a paracomplex in the following sense \cite{Voigtepch}. 
\begin{definition} 
A paracomplex $ C = C_0 \oplus C_1 $ is a given by $ \ayd $-modules $ C_0 $ and $ C_1 $ 
in $ \mathcal{C} $ together with $ \ayd $-maps $ \partial_0: C_0 \rightarrow C_1 $ and $ \partial_1: C_1 \rightarrow C_0 $ such that
$$ 
\partial^2 = \id - T.
$$
A chain map $ \phi: C \rightarrow D $ between two paracomplexes is an $ \ayd $-map from $ C $ to $ D $ that commutes with
the differentials.
\end{definition}
In the same way as for ordinary chain complexes one defines homotopies, mapping cones and suspensions. Together with this additional structure 
the homotopy category of paracomplexes becomes a triangulated category. \\
In order to obtain the paracomplexes needed for the definition of equivariant cyclic homology one has to insert certain 
tensor algebras into the equivariant $ X $-complex. Consider the space $ \Omega(A) $ of ordinary differential forms over a $ G $-algebra $ A $. 
By definition, the analytic bornology on $ \Omega(A) $ is the bornology generated by the sets 
\begin{equation*}
S \cup \bigcup_{n = 1}^\infty S(dS)^n \cup (dS)^n
\end{equation*}
for all small sets $ S \subset A $. Unless explicitly stated otherwise, we will always equip $ \Omega(A) $ with the analytic bornology. 
The space $ \Omega(A) $ becomes a separated $ G $-algebra with the Fedosov 
product 
$$
\omega \circ \eta = \omega \eta - (-1)^{|\omega|} d\omega d\eta
$$ 
for homogenous forms $ \omega $ and $ \eta $ and the diagonal action of $ G $. 
By definition, the analytic tensor algebra $ \mathcal{T}A $ of $ A $ is 
the subalgebra of differential forms of even degree. The basic ingredient in the definition 
of periodic, analytic and local cyclic homology is the equivariant $ X $-complex of 
the tensor algebra. \\
To obtain interesting homology groups one has to complete this paracomplex in some way. 
The most evident way of completion, namely the completion of the underlying bornological vector space, 
leads to the analytic theory. Before giving the general definition of equivariant analytic cyclic homology let 
us first discuss the completions needed to define the local and the periodic theory, respectively. \\
For the local theory the ordinary completion is replaced by the derived completion. The derived completion 
is the left derived functor of the completion functor with respect to the localization of 
the homotopy category of paracomplexes at the class of locally contractible paracomplexes. A paracomplex $ C $ is called 
locally contractible provided all chain maps $ P \rightarrow C $ are homotopic to zero whenever $ P $ is a primitive paracomplex. By definition, 
a paracomplex $ P $ is called primitive if its underlying $ \ayd $-module is primitive. 
The class of locally contractible paracomplexes forms a null system in the homotopy category of paracomplexes and 
the corresponding localization is called the local derived category. 
The derived completion of a paracomplex $ C $ is given by 
$$
C^{\mathbb{L}c} \cong P(C)^c
$$
where $ P(C) \rightarrow C $ is a projective resolution with respect to the class of locally contractible paracomplexes. 
For the explicit construction of a projective resolution functor we refer to \cite{Voigtelch}. 
We remark that there is a natural map $ C^{\mathbb{L}c} \rightarrow C^c $ for every paracomplex $ C $. \\
The completion needed for the periodic theory is obtained using quotients of the analytic 
tensor algebra. More precisely, if $ A $ is a separated $ G $-algebra then the natural homomorphism 
$ \tau_A: \mathcal{T}A \rightarrow A $ given by the projection onto 
differential forms of degree zero fits into an extension 
$$
\xymatrix{
   \mathcal{J}A\; \ar@{>->}[r] & \mathcal{T}A \ar@{->>}[r] & A \\
 }
$$
of separated $ G $-algebras where $ \mathcal{J}A $ denotes the kernel of the map $ \tau_A $. The $ n $-th power $ (\mathcal{J}A)^n $ of 
the ideal $ \mathcal{J}A $ can 
be identified with the space of even differential forms of degree greater or equal $ 2n $. By definition, the periodic 
tensor algebra of $ A $ is the projective system of $ G $-algebras
$$ 
\mathcal{T}A/(\mathcal{J}A)^\infty = (\mathcal{T}A/(\mathcal{J}A)^n)_{n \in \mathbb{N}} 
$$
which can be viewed as a pro-$ G $-algebra in a natural way. Accordingly, the equivariant $ X $-complex of the periodic tensor algebra is 
a pro-paracomplex. In the context of pro-modules and pro-paracomplexes, morphisms are always understood in the sense of pro-categories. \\
We shall now give the definition of equivariant periodic, analytic and local cyclic homology. 
\begin{definition} \label{defanalytic} 
Let $ G $ be a totally disconnected group and let $ A $ and $ B $ be separated $ G $-algebras. 
The equivariant analytic cyclic homology of $ A $ and $ B $ is
\begin{equation*}
HA^G_*(A,B) =
H_*(\SHom_G(X_G(\mathcal{T}(A \otimes \mathcal{K}_G))^c, 
X_G(\mathcal{T}(B \otimes \mathcal{K}_G))^c)).
\end{equation*}
The equivariant local cyclic homology $ HL^G_*(A,B) $ is obtained by replacing the completion with the left derived completion. 
For equivariant periodic cyclic homology $ HP^G_*(A,B) $ one has to replace instead the analytic tensor 
algebras by periodic tensor algebras. 
\end{definition}
Here $ \mathcal{K}_G $ is a certain subalgebra of the algebra of compact operators on $ L^2(G) $. By definition, it consists 
of integral operators with kernels in $ \D(G \times G) $. We point out that, despite of the fact that both entries are only 
paracomplexes, the $ \Hom $-complexes occuring on the right hand side are complexes in the usual sense. \\
It is sometimes useful to extend the previous definition to bornological algebras $ A $ and $ B $ with an action of $ G $
by bounded automorphisms which is not necessarily smooth.  
This is done simply by replacing $ A $ and $ B $ on the right hand side by their smoothings $ \Smooth(A) $ and $ \Smooth(B) $, 
respectively. \\
For our considerations below it is crucial to compare the paracomplex $ X_G(\mathcal{T}A) $ with 
a certain $ B + b $-complex. More precisely, let $ A $ be a separated $ G $-algebra and consider the bornology on 
$ \Omega_G(A) $ generated by the sets 
\begin{equation*}
D \otimes S \cup D \otimes [S] dS \cup \bigcup_{n = 1}^\infty n!\, D \otimes [S][dS] (dS)^{2n} 
\end{equation*}
for small subsets $ D \subset \mathcal{O}_G $ and $ S \subset A $. Here we write $ [R] $ for the set $ R \cup \{1 \}$ 
where $ 1 $ is the unit in the unitarization. 
This bornology is called the transposed analytic bornology. 
We denote by $ \Omega^{\an}_G(A) $ the space $ \Omega_G(A) $ equipped with the transposed analytic bornology. 
This space is an $ \ayd $-module in a natural way and there is a bornological isomorphism 
$ X_G(\mathcal{T}A) \cong \Omega^{\an}_G(A) $ for every $ G $-algebra $ A $. 
Moreover the operators $ b $ and $ B $ induce bounded maps on $ \Omega^{\an}_G(A) $. In particular, together with the boundary 
$ B + b $ the space $ \Omega^{\an}_G(A) $ becomes a paracomplex. The following theorem is proved in \cite{Voigtelch}. 
\begin{theorem}\label{homotopyeq}
For every separated $ G $-algebra $ A $ there exists a natural 
homotopy equivalence between the paracomplexes $ X_G(\mathcal{T}A) $ and $ \Omega^{\emph{\an}}_G(A) $. 
\end{theorem}
As a consequence, in the definition of the analytic and local theory one may replace 
the equivariant $ X $-complexes by the corresponding $ B + b $-complexes. An analogous result holds for the periodic theory \cite{Voigtepch}.\\
Finally, we remark that according to the decomposition of an $ \ayd $-module into an elliptic and a hyperbolic part discussed in section \ref{secsmoothcov} 
we obtain corresponding decompositions of the paracomplexes defining the different cyclic theories. This will be used 
in section \ref{secindcyclic} when we study induction in equivariant cyclic homology. 

\section{Restriction and induction of modules and algebras} \label{secindmodules}

In this section we discuss restriction and compact induction of $ G $-modules, $ G $-algebras and $ \ayd $-modules. 
We only consider these constructions in the particular cases which are needed in the sequel. For a detailed account to 
induction and restriction of smooth representations of locally compact groups we refer to \cite{Meyersmoothrep}. 
We remark that in the context of equivariant sheaves for totally disconnected groups these constructions are treated in \cite{Schneider}. \\
We begin with compact induction of smooth representations. Let $ K $ be a compact open subgroup of a totally disconnected group $ G $
and let $ V $ be a separated $ K $-module. The compactly induced module $ \ind_K^G(V) $ of $ V $ is given by 
$$
\ind_K^G(V) = \{f \in \D(G, V)| \, f(ts) = s^{-1} \cdot f(t)\; \text{for all} \; t \in G, s \in K \}
$$
where $ G $ acts by left translations. Since $ K $ is open there exists a natural $ K $-equivariant map 
$ \theta_V: V \rightarrow \ind_K^G(V) $ given by $ \theta_V(v)(t) = \chi(t) t^{-1} \cdot v $ where $ \chi $ is the characteristic 
function of the set $ K $. \\
One may equivalently describe the compactly induced module as follows. Since $ K $ is an open subgroup of $ G $ there 
is an evident map $ \D(K) \rightarrow \D(G) $ given by extending functions by zero. We may thus 
view $ \D(K) $ as a subalgebra of $ D(G) $ if the Haar measure on $ K $ is chosen to be the restriction of 
the Haar measure on $ G $. Then there is a natural isomorphism $ \Phi: \D(G) \otimes_{\D(K)} V \rightarrow \ind_K^G(V) $ 
given by 
$$
\Phi(f)(t) = \int_K s \cdot f(ts) ds
$$
for the class of $ f \in \D(G,V) = \D(G) \otimes V $. \\
If $ W $ is a separated $ G $-module then restriction of the action to $ K $ yields the structure of a $ K $-module on $ W $. 
We write $ \res^G_K(W) $ for the $ K $-module obtained this way. Compact induction and restriction are related 
as follows. 
\begin{prop}\label{indmodules}
Let $ K $ be a compact open subgroup of the totally disconnected group $ G $. Then the map $ \theta_V $ induces a natural isomorphism
$$
\Hom_G(\ind_K^G(V),W) \cong \Hom_K(V,\res^G_K(W)) 
$$
for all separated $ K $-modules $ V $ and separated $ G $-modules $ W $. 
\end{prop}
\proof The map $ \Psi: \Hom_K(V,\res^G_K(W)) \rightarrow \Hom_G(\ind_K^G(V),W) $ defined by 
$$
\Psi(\phi)(f) = \sum_{t \in G/K} t \cdot \phi(f(t))
$$ 
is inverse to the canonical map $ \Hom_G(\ind_K^G(V),W) \rightarrow \Hom_K(V,\res^G_K(W)) $ induced by $ \theta_V $. \qed \\
Remark that induction is compatible with completion in the sense that the natural map 
$ \ind_K^G(V)^c \rightarrow \ind_K^G(V^c) $ is an isomorphism for every separated $ K $-module $ V $. \\
Next we consider induction and restriction for algebras. Let again $ K $ be a compact open subgroup of $ G $ and let 
$ B $ be a separated $ K $-algebra. By definition, the induced $ G $-algebra is $ \ind_K^G(B) $ with pointwise multiplication. 
The natural map $ \theta_B: B \rightarrow \ind_K^G(B) $ is a $ K $-equivariant algebra homomorphism. 
Moreover there is a natural $ K $-equivariant homomorphism 
$ \pi_B: \ind_K^G(B) \rightarrow B $ given by $ \pi_B(f) = f(e) $ and satisfying $ \pi_B \theta_B = \id $. 
On the other hand every separated $ G $-algebra $ A $ 
is a $ K $-algebra in a natural way by restriction of the action. We will frequently denote
this $ K $-algebra again by $ A $ instead of $ \res^G_K(A) $. \\
Let us now discuss induction and restriction of $ \ayd $-modules. If $ K $ is a compact open subgroup of 
$ G $ there is a natural $ K $-equivariant inclusion homomorphism  
$ \mathcal{O}_K \rightarrow \mathcal{O}_G $. Moreover we have an algebra homomorphism $ \D(K) \rightarrow \D(G) $ as explained 
above. These maps yield an algebra homomorphism $ \iota: \AYD(K) \rightarrow \AYD(G) $. The induced module $ \ind_K^G(M) $ of a 
separated $ K $-$ \ayd $-module $ M $ is the $ G $-$ \ayd $-module defined by 
$$
\ind_K^G(M) = \AYD(G) \otimes_{\AYD(K)} M 
$$
with the obvious left $ \AYD(G) $-module structure. 
Conversely, let $ M $ be a separated $ G $-$ \ayd $-module. We define the restriction $ \res^G_K(M) $ of $ M $ to $ K $ by 
$$
\res^G_K(M) = \mathcal{O}_K \otimes_{\mathcal{O}_G} M 
$$
where $ \mathcal{O}_K $ is viewed as an $ \mathcal{O}_G $-module by restriction of functions. 
The space $ \res^G_K(M) $ becomes a $ K $-module using the diagonal action. Moreover it is a nondegenerate 
$ \mathcal{O}_K $-module by multiplication of functions in the first tensor factor. In this way $ \res^G_K(M) $ becomes a 
$ K $-$ \ayd $-module. \\
The following proposition describes the relation between induction and restriction of $ \ayd $-modules. 
\begin{prop}\label{covinduction}
Let $ G $ be a totally disconnected group and let $ K \subset G $ be a compact open subgroup. Then there exists a natural isomorphism 
$$
\SHom_G(\ind_K^G(M),N) \cong \SHom_K(M,\res_K^G(N))
$$
for all separated $ K $-$ \ayd $-modules $ M $ and $ G $-$ \ayd $-modules $ N $.
\end{prop}
\proof We define a map $ \alpha: \SHom_K(M,\res_K^G(N)) \rightarrow \SHom_G(\ind_K^G(M),N) $ by
$$
\alpha(\phi)(f \otimes m) = f \cdot \mu(\phi(m))
$$
where $ \mu: \mathcal{O}_K \otimes_{\mathcal{O}_G} N \rightarrow N $ is the multiplication map. 
Conversely, we construct a map $ \beta: \SHom_G(\ind_K^G(M),N) \rightarrow \SHom_K(M,\res_K^G(N)) $ as 
follows. Given a $ G $-$ \ayd $-map $ \phi: \ind_K^G(M) \rightarrow N $ we let $ \beta(\phi) $ 
be the composition 
\begin{equation*} 
M \;\cong 
 \xymatrix{
     \AYD(K) \otimes_{\AYD(K)} M \; \ar@{->}[r]^{\iota \otimes \id } &
         \AYD(G) \otimes_{\AYD(K)} M \ar@{->}[r]^{\qquad \quad \phi} & N \ar@{->}[r] & \res^G_K(N)
     }
\end{equation*}
where the last maps sends $ n \in N $ to $ 1 \otimes n \in \mathcal{O}_K \otimes_{\mathcal{O}_G} N = \res^G_K(N) $.  
For the first isomorphism we use the fact that $ M $ is a nondegenerate $ \AYD(K) $-module. Explicitly, we 
have  
$$
\beta(\phi)(f \cdot m) = 1 \otimes \phi(\iota(f) \otimes m)
$$
for $ f \in\AYD(K) $ and $ m \in M $. 
It is straightforward to check that $ \alpha $ and $ \beta $ are inverse isomorphisms. \qed \\
We will mainly have to work with $ G $-$ \ayd $-modules of the form $ M = \mathcal{O}_G \otimes V $ for some separated 
$ G $-module $ V $. In this case the reduced module $ \res_K^G(M) $ is equal to $ \mathcal{O}_K \otimes \res_K^G(V) $. 
In particular, we may identify $ \res_K^G(X_G(A)) $ with $ X_K(\res^G_K(A)) $ for every separated $ G $-algebra $ A $. \\
Similarly, there is a useful description of the induced modules for $ K $-$ \ayd $-modules of the form $ \mathcal{O}_K \otimes V $ 
for a separated $ K $-module $ V $. Let $ \chi \in \mathcal{O}_G $ be the characteristic function of the set $ K $. Note that every element 
in $ \mathcal{O}_G $ determines a multiplier of $ \AYD(G) $ in a natural way.  
\begin{lemma}\label{indlemma2}
Let $ G $ be a totally disconnected group and let $ K \subset G $ be a compact open subgroup. For every separated $ K $-module $ V $ 
we have a natural isomorphism 
$$
\ind_K^G(\mathcal{O}_K \otimes V) \cong I \otimes_{\D(K)} V 
$$
where $ I \subset \AYD(G) $ denotes the left ideal generated by the multiplier $ \chi $. 
\end{lemma}
\proof First observe that the ideal $ I $ is indeed right $ K $-invariant. 
We define a map $ \alpha: \ind_K^G(\mathcal{O}_K \otimes V) \rightarrow I \otimes_{\D(K)} V $ by 
$$
\alpha(f \otimes h \otimes v) = f \cdot h \otimes v
$$
for $ f \in \AYD(G) $ and $ h \otimes v \in \mathcal{O}_K \otimes V $. 
Conversely, let $ \beta: I \otimes_{\D(K)} V \rightarrow \ind_K^G(\mathcal{O}_K \otimes V) $ be given 
by 
$$
\beta(f \otimes v) = f \otimes \chi \otimes v. 
$$
It is easy to check that $ \alpha $ and $ \beta $ define inverse isomorphisms. \qed \\
Finally, remark that if $ M $ is a paracomplex of $ K $-$ \ayd $-modules then the induced module $ \ind_K^G(M) $ 
is a paracomplex of $ G $-$ \ayd $-modules in a natural way. Similarly, if $ N $ is a paracomplex 
of $ G $-$ \ayd $-modules then its restriction $ \res^G_K(N) $ is a paracomplex of $ K $-$ \ayd $-modules. 
Observe moreover that the isomorphism obtained in proposition \ref{covinduction} is a chain map in this case. 

\section{Induction in equivariant cyclic homology} \label{secindcyclic}

In this section we study the compatibility of equivariant cyclic homology with induction. We 
will restrict attention to the situation which arises in the construction of 
the equivariant Chern character later on. \\
The main result of this section is the following theorem. 
\begin{theorem}\label{homspacetheorem}
Let $ K $ be a compact open subgroup of $ G $. For every separated $ K $-algebra $ A $ and every separated 
$ G $-algebra $ B $ there is a natural homotopy equivalence 
$$
\ind_K^G(X_K(\mathcal{T}(A \otimes \res^G_K(B))))^{\mathbb{L}c} \simeq X_G(\mathcal{T}(\ind_K^G(A) \otimes B))^{\mathbb{L}c}
$$
of paracomplexes of $ G $-$ \ayd $-modules. There 
are analogous homotopy equivalences if the derived completion is replaced by the ordinary completion or if one considers 
periodic tensor algebras instead. 
\end{theorem}
The proof of theorem \ref{homspacetheorem} is divided into several steps. We will carry out the details for 
the analytic tensor algebra and the derived completion, the other assertions are obtained in a similar way. 
For simplicity we assume that the Haar measure 
on $ G $ is normalized such that the measure of the compact set $ K $ is $ 1 $. \\
Consider the natural $ K $-equivariant homomorphism $ \theta_A: A \rightarrow \ind_K^G(A) $ given by $ \theta_A(a)(t) = \chi(t) t^{-1} \cdot a $ where  
$ \chi $ denotes the characteristic function of the set $ K $. Similarly, we have the $ K $-equivariant homomorphism 
$ \pi_A: \ind_K^G(A) \rightarrow A $ given by $ \pi_A(\alpha) = \alpha(e) $. 
Tensoring the previous maps with the identity on $ B $ 
yields $ K $-equivariant homomorphisms $ \theta: A \otimes B \rightarrow \ind_K^G(A) \otimes B $ and $ \pi: \ind_K^G(A) \otimes B \rightarrow A \otimes B $, 
respectively. In particular, $ \theta $ induces a $ K $-equivariant 
homomorphism $ \Theta: \mathcal{T}(A \otimes B) \rightarrow \mathcal{T}(\ind_K^G(A) \otimes B) $ and an 
$ \AYD(K) $-linear chain map
$ X_K(\mathcal{T}(A \otimes B)) \rightarrow X_G(\mathcal{T}(\ind_K^G(A) \otimes B)) $. According to 
proposition \ref{covinduction} this map corresponds to an $ \AYD(G) $-linear chain map
$$ 
i: \ind_K^G(X_K(\mathcal{T}(A \otimes B))) \rightarrow X_G(\mathcal{T}(\ind_K^G(A) \otimes B))
$$ 
and using lemma \ref{indlemma2} we obtain the explicit formula 
\begin{align*}
i(F \otimes x_0)(s) &= \int_G F(s,t) t \cdot \Theta(x_0) dt \\
i(F \otimes x_0 Dx_1)(s) &= \int_G F(s,t) (t \cdot \Theta(x_0)) D(t \cdot \Theta(x_1)) dt 
\end{align*}
for $ i $ where $ F $ is an element in the left ideal $ I \subset \AYD(G) $ generated by the multiplier $ \chi $ and 
$ x_i \in \mathcal{T}(A \otimes B) $. \\
We construct a chain map $ p $ which is going to be the homotopy inverse of $ i $. Using lemma \ref{indlemma2} we 
define $ p: X_G(\mathcal{T}(\ind_K^G(A) \otimes B)) \rightarrow \ind_K^G (X_K(\mathcal{T}(A \otimes B))) $ by 
\begin{align*}
p(f \otimes \alpha_0)(s,t) &= f(s) (t \cdot \chi)(s) \Pi(t^{-1} \cdot \alpha_0) \\
p(f \otimes \alpha_0 D\alpha_1)(s,t) &= f(s) (t \cdot \chi)(s) \Pi(t^{-1} \cdot \alpha_0) D \Pi(t^{-1} \cdot \alpha_1)
\end{align*}
for $ \alpha_i \in \mathcal{T}(\ind_K^G(A) \otimes B) $ where $ \Pi: \mathcal{T}(\ind_K^G(A) \otimes B)) \rightarrow \mathcal{T}(B \otimes A) $ denotes 
the homomorphism induced by $ \pi: \ind_K^G(A) \otimes B \rightarrow A \otimes B $. Remark that $ \Pi(t \cdot \alpha) $ has compact support as a 
function of $ t $ for every $ \alpha \in \mathcal{T}(\ind_K^G(A) \otimes B) $. 
It is easy to check that $ p $ is an $ \ayd $-map, and by definition $ p $ commutes with the boundary operator $ d $. A
straightforward computation shows that  $ p $ commutes with the boundary $ b $ as well. It follows that $ p $ is a chain map. 
We compute 
\begin{align*}
pi(F \otimes x_0)(s,t) &= \int_G F(s,r)(t \cdot \chi)(s) \Pi(t^{-1}r \cdot \Theta(x_0)) dr \\
&= \int_K F(s, tr)(t \cdot \chi)(s) \Pi(r \cdot \Theta(x_0)) dr 
\end{align*}
where the second equality follows from left invariance of the Haar measure and the fact that $ \Pi(r \cdot \Theta(x_0)) $ vanishes for $ r $ outside 
$ K $. Using that $ \Pi $ is $ K $-equivariant we obtain
\begin{align*}
\int_K F(s, tr)(t \cdot &\chi)(s) \Pi(r \cdot \Theta(x_0)) dr = \int_K F(s, t)(t \cdot \chi)(s) \Pi\Theta(x_0) dr \\
&= F(s, t)(t \cdot \chi)(s) x_0 = F(s, t)x_0
\end{align*} 
taking into account that the Haar measure on $ K $ is normalized and that $ F $ is contained in the ideal $ I $. 
In the same way one calculates 
$$
pi(F \otimes x_0 Dx_1)(s,t) = F(s,t)x_0 Dx_1
$$
in degree one. Hence we have proved the relation $ pi = \id $. \\
Let us now assume that $ A $ and $ B $ are unital and show that $ ip $ is homotopic to the identity. For this it 
is convenient to use theorem \ref{homotopyeq} and to work with equivariant differential forms instead of equivariant $ X $-complexes. 
Rewriting the formula for $ i $ yields
$$
i(F \otimes x_0 dx_1\cdots dx_n)(s) = \int_G F(s,t) (t \cdot \theta(x_0)) d(t \cdot \theta(x_1)) \cdots d(t \cdot \theta(x_n)) dt
$$
on $ \ind_K^G(\Omega^n_K(A \otimes B)) $. For the map $ p $ we obtain
$$
p(f \otimes \alpha_0 d\alpha_1 \cdots d\alpha_n)(s,t) = f(s) (t \cdot \chi)(s)  \pi(t^{-1} \cdot \alpha_0) d \pi(t^{-1} \cdot \alpha_1) \cdots 
d\pi(t^{-1} \cdot \alpha_n) 
$$
on $ \Omega^n_G(\ind_K^G(A) \otimes B) $. 
Given $ \alpha = \beta \otimes b \in \ind_K^G(A) \otimes B $ let us denote by 
$ \alpha(t) $ the element $ \beta(t) \otimes b $ in $ A \otimes B $. 
Moreover, for $ t \in G $ we write $ \alpha[t] $ for the element in $ \ind_K^G(A) \otimes B $ given by $ \alpha[t](s) = \alpha(s) \chi(t^{-1}s) $. 
That is, $ \alpha[t] $ is equal to the function $ \alpha $ on the coset $ t K $ and zero elsewhere. 
Using this notation we obtain
\begin{align*}
ip(f &\otimes \alpha_0 d\alpha_1 \cdots d\alpha_n)(s) \\
&= \int_G f(s) (t \cdot \chi)(s) (t \cdot (\theta \pi(t^{-1} \cdot \alpha_0)) 
d(t \cdot \theta \pi(t^{-1} \cdot \alpha_1)) \cdots d(t \cdot \theta \pi(t^{-1} \cdot \alpha_n)) dt\\
&= \int_G f(s) (t \cdot \chi)(s)\; \alpha_0[t]
d \alpha_1[t] \cdots d \alpha_n[t] dt. 
\end{align*}
This implies
\begin{align*}
(\id - ip)&(f \otimes \alpha_0 d\alpha_1 \cdots d\alpha_n)(s) = \int (1 - \chi(t_1) \cdots \chi(t_n) (t_0 \cdot \chi)(s)) f(s) \\
& \alpha_0[t_0] d \alpha_1[t_0t_1] \cdots d \alpha_n[t_0 \cdots t_n]\; dt_0 dt_1 \cdots dt_n
\end{align*}
and 
\begin{align*}
(\id - ip)&(f \otimes d\alpha_1 \cdots d\alpha_n)(s) = \int (1 - \chi(t_2) \cdots \chi(t_n) (t_1 \cdot \chi)(s)) f(s)  \\
& d \alpha_1[t_1] \cdots d \alpha_n[t_1 \cdots t_n]\; dt_1 \cdots dt_n.
\end{align*}
Let us abbreviate $ \C = \Omega_G(\ind_K^G(A) \otimes B) $ and define a map $ h: \C \rightarrow \C $ by  
\begin{align*}
h&(f \otimes \alpha_0 d\alpha_1 \cdots d\alpha_n)(s) = \sum_{j = 0}^{n - 1} (-1)^j \int \chi(t_1) \cdots \chi(t_j) \eta(t_{j + 1})\, f(s) \,
\alpha_0[t_0] d \alpha_1[t_0 t_1] \cdots \\
&\qquad \cdots d\alpha_j[t_0 \cdots t_j] d1[t_0 \cdots t_j] d \alpha_{j + 1}[t_0 \cdots t_{j + 1}] \cdots d\alpha_n[t_0 \cdots t_n]\, dt_0 \cdots dt_n \\
&+ (-1)^n \int \chi(t_1) \cdots \chi(t_n)\, f(s) (t_0 \cdot\eta)(s) 
\,\alpha_0[t_0] d\alpha_1[t_0t_1] \cdots \\
& \qquad \cdots d\alpha_n[t_0 \cdots t_n] d1[t_0 \cdots t_n] \; dt_0 \cdots dt_n 
\end{align*}
where $ \eta(t) = 1 - \chi(t) $ as well as 
\begin{align*}
h(&f \otimes d\alpha_1 \cdots d\alpha_n)(s) = \sum_{j = 1}^{n - 1} (-1)^j \int \chi(t_2) \cdots \chi(t_j) \eta(t_{j + 1}) \, f(s) \,
d \alpha_1[t_1] d\alpha_2[t_1t_2]\cdots \\
&\qquad \cdots d\alpha_j[t_1 \cdots t_j] d1[t_1 \cdots t_j] d\alpha_{j + 1}[t_1 \cdots t_{j + 1}] 
\cdots d\alpha_n[t_1 \cdots t_n] \; dt_1 \cdots dt_n \allowdisplaybreaks[3] \\
&+ \sum_{j = 1}^{n - 1} (-1)^{j + n} \int \chi(t_2) \cdots \chi(t_j) \eta(t_{j + 1}) \, f(s) \eta(t_1^{-1}s^{-1} t_1 \cdots t_n) \\
&\qquad s^{-1} \cdot \alpha_n[t_1 \cdots t_n] d s^{-1} \cdot 1[t_1 \cdots t_n] d \alpha_1[t_1] \cdots \\
&\qquad \cdots d\alpha_j[t_1 \cdots t_j] d1[t_1 \cdots t_j] d\alpha_{j + 1}[t_1 \cdots t_{j + 1}] \cdots d\alpha_{n - 1}[t_1 \cdots t_{n - 1}]\; dt_1 \cdots dt_n \allowdisplaybreaks[2] \\
& + (-1)^n \int \chi(t_2) \cdots \chi(t_n)\, f(s) (t_1 \cdot \eta)(s) 
\, d\alpha_1[t_1] d \alpha_2[t_1t_2] \cdots \\
& \qquad \cdots d\alpha_n[t_1 \cdots t_n] d1[t_1 \cdots t_n] \; dt_1 \cdots dt_n \allowdisplaybreaks[2] \\
& + \int \chi(t_2) \cdots \chi(t_n)\, f(s) (t_1 \cdot\eta)(s) 
\,d s^{-1} \cdot 1[t_1] d\alpha_1[t_1] \cdots d \alpha_n[t_1 \cdots t_n] \; dt_1 \cdots dt_n.
\end{align*}
It is easy to check that $ h $ is an $ \ayd $-map and bounded with respect to the analytic bornology. A lengthy but straightforward 
computation yields the relation 
$$
bh + hb = \id - ip
$$
on $ \C $ which shows that $ ip $ is homotopic to the identity with respect to the Hochschild operator. 
Moreover we have $ h i = 0 $ and $ p h = 0 $. \\
Recall that $ M \otimes_T \mathbb{C} $ denotes the space of $ T $-coinvariants of an $ \ayd $-module $ M $. 
We are interested in the space $ \C \otimes_T \mathbb{C} $ of coinvariants of $ \C $. 
The transposed analytic bornology on this space is generated by sets of the form 
\begin{equation*}
D \otimes S \cup D \otimes [S] dS \cup \bigcup_{n = 1}^\infty n!\, D \otimes [S][dS] (dS)^{2n} 
\end{equation*}
where $ D \subset \mathcal{O}_G $ and $ S \subset \ind_K^G(A) \otimes B $ are 
small. Actually it suffices to consider sets $ D $ and $ S $ of a special form. 
Let us fix a $ K $-invariant small disk $ R $ of $ A \otimes B $ containing the 
identity element. We may assume that $ S $ consists of all $ \alpha $ in $ \ind_K^G(A) \otimes B $ 
with support in a $ K $-invariant compact subset $ L $ of $ G $ such that 
the element $ \alpha(t) $ is contained in $ R $ for all $ t \in G $. 
For the set $ D $ we may assume that it is closed under multiplication by all functions 
$ t \cdot \eta $ for $ t \in L $. Let us call the pair $ D $ and $ S $ special if these conditions 
are satisfied. \\
We define endomorphisms $ T_j $ of $ \Omega^n_G(\ind_K^G(A) \otimes B) \otimes_T \mathbb{C} $ by the
formula 
\begin{equation*}
T_j(f \otimes \alpha_0 d\alpha_1 \cdots d\alpha_n)(s)
= f(s)\, (s^{-1} \cdot \alpha_0)d(s^{-1} \cdot \alpha_1) \cdots d(s^{-1} \cdot \alpha_{j - 1}) d\alpha_j \cdots d\alpha_n,
\end{equation*}
and for small subsets $ D $ and $ S $ as before we set
\begin{align*}
D\Omega^n(S) = \bigcup_{j = 0}^n T_j(D \otimes [S] (dS)^n) 
\end{align*}
for $ n > 0 $ as well as $ D\Omega^0(S) = D \otimes S $. Using this notation we obtain the following estimate. 
\begin{lemma}\label{homotopybounded1}
Let $ D $ and $ S $ be a pair of special small sets. There exists a positive constant $ m $ dependent on $ S $ such that 
$ (hB)^j(x) $ can be written as a sum of at most 
$$
m^{n + 1} (n + 1) (n + 3) \cdots (n + 2j - 1) 2
$$
terms in $ D\Omega^{n + 2j}(S) $ for any $ x \in D\Omega^n(S) $ and $ j \in \mathbb{N} $. 
\end{lemma}
\proof Let us fix an element $ x \in D\Omega^n(S) $. Then, by definition, the element $ B(x) $ is given as a sum of 
$ n + 1 $ terms contained in $ D\Omega^{n + 1}(S) $. 
Similarly, $ h(x) $ is given as a finite sum of elementary terms. 
Let $ m $ be the number of elements in the image of $ L $ in $ G/K $ under the canonical projection where $ L \subset G $ is the compact subset 
occuring in the definition of $ S $. Using that $ D $ and $ S $ 
are special it is easily seen that each 
term in the formula for $ h(x) $ can be written as a sum of at most
$ m^{n + 1} $ elements belonging to $ D\Omega^{n + 1}(S) $. 
It suffices to estimate the number of nonzero summands in $ (hB)^j(x) $ for all $ j $. 
In the definition of $ h $ for elements in the image of the map $ d $ only the 
terms in the second sum are not contained in the kernel of $ d $. It follows that $ (hB)^j(x) $ 
can be written as a sum of at most $ m^{n + 1} (n + 1) (n + 3) \cdots (n + 2j - 1) 2 $ terms 
in $ D\Omega^{n + 2j}(S) $ as claimed. \qed \\
With the notation as above we let $ D\Omega(S) $ be the circled convex hull of the set
$$
\bigcup_{n = 0}^\infty n!\, D\Omega^{2n}(S) \cup \bigcup_{n = 0}^\infty n!\, D\Omega^{2n + 1}(S)
$$
in $ \C \otimes_T \mathbb{C} $. If $ D $ and $ S $ are small then the set $ D\Omega(S) $ is small in the transposed analytic bornology. 
Using lemma \ref{homotopybounded1} we construct a map on the coinvariants of the derived 
completion $ \C^{\mathbb{L}c} $ as follows. 
\begin{prop}\label{homotopybounded2}
On the space of coinvariants $ \C^{\mathbb{L}c} \otimes_T \mathbb{C} $ of the derived completion of the paracomplex $ \C $ the series
\begin{equation*} 
\sum_{j = 0}^\infty (-1)^j (hB)^j
\end{equation*} 
determines an endomorphism of $ \ayd $-modules. 
\end{prop}
\proof Every small subset of $ \C \otimes_T \mathbb{C} $ is contained in a 
set of the form $ D\Omega(S) $ where $ D $ and $ S $ are a special pair of small sets in the sense above. 
Using lemma \ref{homotopybounded1} we see that the partial sums in the above series
satisfy the Cauchy criterion if they are viewed as bounded linear maps from the normed space 
$ \bra D\Omega(S) \ket $ into the normed space $ \bra D\Omega(\lambda S) \ket $ provided $ \lambda $ 
is chosen such that $ \lambda^{n + 2j} \geq 2 m^{n + 1} 4^j $ for all $ n $ and $ j $. 
It follows that the series converges to a bounded $ \ayd $-map from 
$ \AYD(G) \bra D\Omega(S) \ket $ into the completion of $ \AYD(G) \bra D\Omega(\lambda S) \ket $. 
In this way one obtains an endomorphism of $ (\C \otimes_T \mathbb{C})^{\mathbb{L}c} $. 
Inspecting the construction of the derived completion yields a natural isomorphism $ \C^{\mathbb{L}c} \otimes_T \mathbb{C} \cong 
(\C \otimes_T \mathbb{C})^{\mathbb{L}c} $. This proves the assertion. \qed  \\
Let us write $ C = \C^{\mathbb{L}c} $ and consider the map $ k: C \otimes_T \mathbb{C} \rightarrow C \otimes_T \mathbb{C} $ 
obtained in proposition \ref{homotopybounded2}. Using this map we are able to apply 
the perturbation lemma \cite{Kassel}. More precisely, together with the previously established properties of 
the maps $ i, p $ and $ h $ we obtain the relation 
$$ 
ip - \id = [\partial, k]
$$ 
on $ C \otimes_T \mathbb{C} $ where $ \partial $ denotes the boundary map $ B + b $. 
For the precise formulation and a proof of the perturbation lemma in the context of paramixed complexes we refer to \cite{Voigtepch}. \\
Let us denote the natural projection $ C \rightarrow C \otimes_T \mathbb{C} $ by $ q $.
It is not hard to check that there exists an $ \ayd $-map $ K: C \rightarrow C $ lifting $ k $. As a consequence
we see in particular that the image of 
\begin{equation*}
\delta = \id - ip + [\partial, K]
\end{equation*}
is contained in the kernel $ (\id - T)C $ of $ q $. 
\begin{lemma}\label{deltazero}
The chain map $ \delta: C \rightarrow C $ is homotopic to zero. Consequently, the chain map 
$ ip: C \rightarrow C $ is homotopic to the identity. 
\end{lemma}
\proof We have to consider the elliptic and hyperbolic part of $ C $ separately. On the elliptic part we have 
$ (\id - T)C_{\elli} = (\id - E)C_{\elli} $ where $ E $ is the canonical projection on the $ T $-invariant part. 
The argument given in the proof of proposition 3.4 in 
\cite{Voigtbs} shows that the natural inclusion $ (\id - E)C_{\elli} \rightarrow C_{\elli} $ is homotopic to zero. 
Hence the chain map $ \delta $ is homotopic to zero on the elliptic part. \\
Now consider the hyperbolic part. 
According to proposition \ref{indmodules} the $ G $-module $ \ind_K^G(A) $ 
is projective. Hence the $ G $-module $ \ind_K^G(A) \otimes B $ is projective as well, and it follows that 
$ \C $ is a projective $ \ayd $-module \cite{Voigtepch}. It is then easy to check that
$$ 
(\id - T): C_{\hyp} \rightarrow (\id - T)C_{\hyp} 
$$ 
is an isomorphism. Hence there exists a chain map $ \Delta: C_{\hyp} 
\rightarrow C_{\hyp} $ such that $ \delta = (\id - T)\Delta $ on $ C_{\hyp} $. 
Since $ \id - T $ is homotopic to zero this shows that the restriction 
of $ \delta $ to the hyperbolic part is homotopic to zero. 
We conclude that $ ip $ is homotopic to the identity on $ C $. \qed \\
Lemma \ref{deltazero} completes the proof of theorem \ref{homspacetheorem} in the case that the algebras $ A $ and $ B $ are unital. 
Now let $ A $ and $ B $ be arbitrary and consider the canonical split extension $ 0 \rightarrow A \rightarrow A^+ \rightarrow \mathbb{C} \rightarrow 0 $ 
of $ K $-algebras. The corresponding extension 
$$
0 \rightarrow \ind_K^G(A) \rightarrow \ind_K^G(A^+) \rightarrow \ind_K^G(\mathbb{C}) \rightarrow 0 
$$
of $ G $-algebras has an equivariant splitting as well, and similarly we have the split
extension $ 0 \rightarrow B \rightarrow B^+ \rightarrow \mathbb{C} \rightarrow 0 $ 
of $ G $-algebras. Using the result in the unital case we obtain the assertion for general $ A $ and $ B $ by considering 
tensor products of the algebras in these extensions and 
excision for $ K $-equivariant and $ G $-equivariant local cyclic homology \cite{Voigtelch}.
This finishes the proof of theorem \ref{homspacetheorem}. \\
Let us discuss some consequences of theorem \ref{homspacetheorem}. First of all, we deduce that induction from compact 
open subgroups descends to the level of equivariant cyclic homology.
\begin{theorem}\label{cyclicinduction}
Let $ K \subset G $ be a compact open subgroup. Then there exists a transformation 
$$
\ind^G_K: HL^K_*(A,B) \rightarrow HL^G_*(\ind^G_K(A), \ind^G_K(B))
$$
which is compatible with the composition product. Analogous transformations exist for the analytic and periodic 
theories. 
\end{theorem}
\proof Using lemma \ref{indlemma2} one checks that induction is compatible with derived 
completion in the sense that there is a natural isomorphism 
$$
\ind_K^G(X_K(\mathcal{T}(C \otimes \mathcal{K}_G))^{\mathbb{L}c}) \cong \ind_K^G(X_K(\mathcal{T}(C \otimes \mathcal{K}_G)))^{\mathbb{L}c}
$$
for every $ K $-algebra $ C $. 
Hence theorem \ref{homspacetheorem} yields a natural homotopy equivalence 
$$
\ind_K^G(X_K(\mathcal{T}(C \otimes \mathcal{K}_K))^{\mathbb{L}c}) \simeq X_G(\mathcal{T}(\ind_K^G(C) \otimes \mathcal{K}_G))^{\mathbb{L}c}
$$
where we take into account stability of $ K $-equivariant cyclic homology \cite{Voigtelch}.
As a consequence one obtains easily the desired induction homomorphism for the local theory. 
The assertions for the analytic and the periodic theories are proved in the same way. \qed \\
Furthermore we obtain the following induction isomorphism. 
\begin{theorem}
Let $ K \subset G $ be a compact open subgroup. For every separated $ K $-algebra $ A $ and every separated $ G $-algebra $ B $ there 
exists a natural isomorphism
$$
HL^G_*(\ind^G_K(A), B) \cong HL^K_*(A, \res^G_K(B)).
$$
Analogous isomorphisms hold for the analytic and periodic theories. 
\end{theorem}
\proof Using the fact that induction is compatible with derived completion this follows again by the homotopy equivalence 
$$
\ind_K^G(X_K(\mathcal{T}(A \otimes \mathcal{K}_G)))^{\mathbb{L}c} \simeq X_G(\mathcal{T}(\ind_K^G(A) \otimes \mathcal{K}_G))^{\mathbb{L}c}
$$
obtained in theorem \ref{homspacetheorem} together with proposition \ref{covinduction} and stability. In the same way one obtains the assertions 
for the analytic and the periodic theories. \qed 

\section{The equivariant Chern character} \label{secchern}

In this section we construct the equivariant Chern character with values in bivariant equivariant cohomology 
in the sense of Baum and Schneider. \\
The first ingredient in this construction is the equivariant Chern-Connes character into
equivariant local cyclic homology obtained in \cite{Voigtelch}. 
\begin{theorem}\label{localchern}
Let $ A $ and $ B $ be separable $ G $-$ C^* $-algebras. Then there exists a transformation  
\begin{equation*}
\ch^G_*: KK^G_*(A,B) \rightarrow HL^G_*(A, B) 
\end{equation*}
which is compatible with the Kasparov product in $ KK^G_* $ and the composition product 
in $ HL^G_* $. This transformation maps elements in $ KK^G_0(A,B) $ induced by equivariant 
$ * $-homomorphisms from $ A $ to $ B $ to the corresponding elements in $ HL^G_0(A,B) $. 
\end{theorem}
We remark that in the context of local cyclic homology all $ G $-$ C^* $-algebras are equipped with the precompact bornology. 
Moreover, it is understood that the smoothing functor $ \Smooth $ is applied to the resulting 
bornological algebras since, unless $ G $ is discrete, the group action on a $ G $-$ C^* $-algebra is not smooth in general. 
For more information we refer to \cite{Voigtelch}. \\
We shall apply theorem \ref{localchern} in the case that $ A $ and $ B $ are algebras of functions on
$ G $-simplicial complexes in the sense of \cite{Voigtbs}. By definition, a $ G $-simplicial complex is a simplicial complex equipped with smooth, 
simplicial and type-preserving action of the group $ G $. In the sequel we will only consider $ G $-simplicial complexes that have at most 
countably many simplices. In order to determine the right hand side 
of the Chern-Connes character further in this situation we have to pass to smooth functions. A suitable notion of smooth functions 
on simplicial complexes was introduced in \cite{Voigtbs}. For a locally finite $ G $-simplicial complex $ Z $ the resulting 
algebra $ C^\infty_c(Z) $ of regular smooth functions with compact support is a $ G $-algebra in a natural way. 
Moreover, it is shown in \cite{Voigtelch} that the inclusion $ C^\infty_c(Z) \rightarrow \Smooth(C_0(Z)) $ induces 
an invertible element in the local theory provided $ Z $ is finite dimensional. It follows that there is an isomorphism
\begin{equation*}
HL^G_*(C_0(X), C_0(Y)) \cong HL^G_*(C^\infty_c(X), C^\infty_c(Y)) 
\end{equation*}
for all finite dimensional locally finite $ G $-simplicial complexes $ X $ and $ Y $ which is compatible with the composition product. \\
The next aim is to determine the equivariant local cyclic homology of algebras of regular smooth functions. 
In order to do this we shall use two auxiliary homology theories. Firstly, we define $ HLA^G_*(A,B) $ by 
$$
HLA^G_*(A,B) = H_*(\SHom_G(X_G(\mathcal{T}(A \otimes \mathcal{K}_G))^{\mathbb{L}c}, 
X_G(\mathcal{T}(B \otimes \mathcal{K}_G))^c))
$$
for all $ G $-algebras $ A $ and $ B $. By construction, this theory is a combination of equivariant local and analytic cyclic homology. 
There is a homomorphism
$$
HL^G_*(A,B) \rightarrow HLA^G_*(A,B)
$$
induced by the canonical map from the derived completion to the ordinary completion in the second variable. 
According to a result from \cite{Voigtelch} the natural chain map $ X_G(\mathcal{T}(B \otimes \mathcal{K}_G))^{\mathbb{L}c} 
\rightarrow X_G(\mathcal{T}(B \otimes \mathcal{K}_G))^c $ is an isomorphism in the local derived category 
provided $ B $ is a Schwartz space satisfying the approximation property. This yields the following assertion. 
\begin{prop}\label{HLAcomp}
Assume that $ B $ is a $ G $-algebra whose underlying bornological vector space is a Schwartz space satisfying the approximation property. Then the natural map 
$$ 
HL^G_*(A,B) \rightarrow HLA^G_*(A,B) 
$$ 
is an isomorphism for every separated $ G $-algebra $ A $. 
\end{prop}
In particular, proposition \ref{HLAcomp} applies in the case that $ B = C^\infty_c(Y) $ for a finite dimensional locally finite $ G $-simplicial complex $ Y $. \\
Let us define another homology theory $ HPA^G_*(A,B) $ by 
$$
HPA^G_*(A,B) = H_*(\SHom_G(X_G(\mathcal{T}(A \otimes \mathcal{K}_G)/\mathcal{J}(A \otimes \mathcal{K}_G)^\infty), 
X_G(\mathcal{T}(B \otimes \mathcal{K}_G))^c))
$$
for all $ G $-algebras $ A $ and $ B $. This theory provides a link between equivariant periodic cyclic homology and 
$ HAL^G_* $. In particular, there is a natural transformation 
$$
HPA^G_*(A,B) \rightarrow HLA^G_*(A,B)
$$
induced by the canonical map from the derived completion to the ordinary completion and the projection from the analytic tensor algebra 
to the periodic tensor algebra in the first variable. 
\begin{prop} \label{theoremlambda}
Let $ X $ be a $ G $-finite proper $ G $-simplicial complex and let $ B $ be any $ G $-algebra. Then the natural map
$$
HPA^G_*(C^\infty_c(X), B) \rightarrow HLA^G_*(C^\infty_c(X), B) 
$$ 
is an isomorphism. 
\end{prop}
\proof Observe that the transformation 
$ HPA^G_*(A, B) \rightarrow HLA^G_*(A, B) $ is compatible with the boundary maps 
in the six-term exact sequences associated to an extension of algebras. Using induction on the 
dimension of $ X $ and excision \cite{Voigtepch}, \cite{Voigtelch} it suffices to prove the assertion in the case that $ X = G/K $ 
is a homogenous space for a compact open subgroup $ K $ of $ G $. Remark that we have $ C^\infty_c(X) = \ind_K^G(\mathbb{C}) $ 
in this situation and consider the commutative diagram 
\begin{equation*}
\xymatrix{
X_G(\mathcal{T}(\ind_K^G(\mathbb{C}) \otimes \mathcal{K}_G))^{\mathbb{L}c} \ar@{->}[r] \ar@{<-}[d]^\simeq & 
X_G(\mathcal{T}(\ind_K^G(\mathbb{C}) \otimes \mathcal{K}_G)/\mathcal{J}(\ind_K^G(\mathbb{C}) \otimes \mathcal{K}_G)^\infty) 
 \ar@{<-}[d]^\simeq  \\
\ind_K^G(X_K(\mathcal{T}(\mathcal{K}_G)))^{\mathbb{L}c} \ar@{->}[r]\ar@{<-}[d]^\simeq & \ind_K^G(X_K(\mathcal{T}(\mathcal{K}_G)/\mathcal{J}(\mathcal{K}_G)^\infty) \ar@{<-}[d]^\simeq \\
\ind_K^G(\mathcal{O}_K[0])^{\mathbb{L}c} \ar@{->}[r]^\cong & \ind_K^G(\mathcal{O}_K[0])
 }
\end{equation*}
where the upper vertical homotopy equivalences are obtained using theorem \ref{homspacetheorem} and the 
lower vertical homotopy equivalences follow from stability in equivariant cyclic homology \cite{Voigtepch}, \cite{Voigtelch}. 
According to the definition of the derived completion the lower horizontal map is an isomorphism. 
Hence the remaining horizontal arrows in this diagram are homotopy equivalences. 
This yields the assertion for $ X = G/K $ and finishes the proof. \qed \\
Applying the natural projection from the analytic to the periodic tensor algebra in the second variable we obtain a transformation 
$$
HPA^G_*(A,B) \rightarrow HP^G_*(A,B)
$$
for all separated $ G $-algebras $ A $ and $ B $. 
\begin{prop} \label{theoremrho}
Let $ X $ and $ Y $ be $ G $-simplicial complexes where $ X $ is $ G $-finite and proper and $ Y $ is finite dimensional and 
locally finite. Then the natural map
$$
HPA^G_*(C^\infty_c(X), C^\infty_c(Y)) \rightarrow HP^G_*(C^\infty_c(X), C^\infty_c(Y))
$$ 
is an isomorphism. 
\end{prop}
\proof Using excision in the first variable we can reduce the assertion to the case that 
$ X = G/K $ is a homogenous space where $ K $ is a compact open subgroup of $ G $. Due to 
theorem \ref{homspacetheorem} it thus suffices to consider the case that $ G $ is compact and that $ X $ is a point. 
Observe that for a compact group $ G $ the chain complexes defining 
$ HPA^G_*(\mathbb{C},B) $ and $ HP^G_*(\mathbb{C}, B) $ for a $ G $-algebra $ B $ are 
obtained by taking the $ G $-invariant part of $ X_G(\mathcal{T}(B))^c $ and 
$ X_G(\mathcal{T}(B)/\mathcal{J}(B)^\infty) $, respectively. \\
We will now use induction on the dimension of $ Y $. If $ \dim(Y) = 0 $ the space $ Y $ is 
the disjoint union of a family of homogenous spaces. Accordingly, we have 
$$
C^\infty_c(Y) = \bigoplus_{j \in J} C(G/K_j) 
$$ 
for some countable index set $ J $ in this case where each $ K_j $ is an open subgroup of the compact group $ G $. 
If the set $ J $ is finite it follows in the same way as in the proof of proposition \ref{theoremlambda}
that the natural map 
$$
X_G(\mathcal{T}(C^\infty(Y)))^c \rightarrow 
X_G(\mathcal{T}(C^\infty(Y))/\mathcal{J}(C^\infty(Y))^\infty)
$$
is a homotopy equivalence. Observe that due to stability we may drop the algebra $ \mathcal{K}_G $ in our 
arguments since $ G $ is compact. 
In order to treat the case of an infinite index set $ J $ we use a direct limit argument. 
More precisely, since completion commutes with direct limits the natural map
$$
\bigoplus_{j \in J} HPA^G_*(\mathbb{C}, C(G/K_j)) \rightarrow HPA^G_*(\mathbb{C}, C^\infty_c(Y)) 
$$
is an isomorphism. For the periodic theory there is an analogous isomorphism according to the results obtained 
in \cite{Voigtbs}. Using the case of finite $ J $ treated before it follows that the claim is true for 
arbitrary $ Y $ of dimension $ 0 $. As in the proof of proposition \ref{theoremlambda} the induction step 
is carried out using excision. \qed \\
The link between equivariant cyclic homology and the theory of Baum and Schneider is provided by the main result of \cite{Voigtbs}. 
\begin{theorem}\label{HPbsiso} Let $ G $ be a totally disconnected group and let $ X $ and $ Y $ be finite dimensional locally 
finite $ G $-simplicial complexes. If the action of $ G $ on $ X $ is proper there exists an isomorphism
\begin{equation*}
HP^G_*(C^\infty_c(X), C^\infty_c(Y)) \cong 
\bigoplus_{j \in \mathbb{Z}} H^{* + 2j}_G(X,Y).
\end{equation*}
\end{theorem}
Combining theorem \ref{localchern} with propositions \ref{HLAcomp}, \ref{theoremlambda}, \ref{theoremrho} and theorem \ref{HPbsiso} we
obtain the following result. 
\begin{theorem}\label{chernbs}
Let $ G $ be a totally disconnected group and let $ X $ and $ Y $ be finite dimensional locally finite $ G $-simplicial 
complexes. If $ X $ is proper and $ G $-finite there exists an equivariant Chern character 
\begin{equation*}
ch_*^G: KK^G_*(C_0(X),C_0(Y)) \rightarrow \bigoplus_{j \in \mathbb{Z}} H^{* + 2j}_G(X,Y)
\end{equation*}
which becomes an isomorphism after tensoring the left hand side with $ \mathbb{C} $. 
\end{theorem}
\proof It is not very hard to verify that the character $ ch^G_* $ is natural with respect to proper equivariant simplicial maps in both variables. 
In order to show that it becomes an isomorphism after tensoring the left hand side with $ \mathbb{C} $ we proceed
in the same way as in the proof of proposition \ref{theoremlambda} and proposition \ref{theoremrho} to reduce to the case that 
$ G $ is compact and $ X $ and $ Y $ are one-point spaces. In this case the assertion follows from proposition 12.4 in \cite{Voigtelch}. \qed \\
Our arguments above show in fact that for $ A = C_0(X) $ and $ B = C_0(Y) $ satisfying 
the assumptions of theorem \ref{chernbs} the Chern-Connes character with values in equivariant 
local cyclic homology is an isomorphism after tensoring the left hand side with $ \mathbb{C} $. Moreover, all the different variants of 
equivariant cyclic homology for the corresponding algebras of smooth functions agree in this case. \\ 
Let us discuss some consequences of theorem \ref{chernbs}. The equivariant $ K $-homology of a proper $ G $-simplicial 
complex $ X $ is defined by 
$$
K^G_*(X) = \varinjlim_{F \subset X} KK^G_*(C_0(F), \mathbb{C})
$$
where the limit is taken over all $ G $-finite subcomplexes $ F $ of $ X $. 
Note that every $ G $-finite subcomplex of $ X $ is locally finite and hence locally compact. 
In a similar way we define the equivariant homology of $ X $ by 
$$
H^G_*(X) = \varinjlim_{F \subset X} HP^G_*(C^\infty_c(F), \mathbb{C}).
$$
if $ X $ is a proper $ G $-simplicial complex. We remark that, taking into account the above considerations, one could as 
well take equivariant analytic or local 
cyclic homology in this definition. Moreover, using theorem \ref{HPbsiso} and the work of Baum and Schneider \cite{BS}
it follows that the groups $ H^G_*(X) $ are naturally isomorphic to the cosheaf homology of $ X $ as defined 
in \cite{BCH}. \\
Taking direct limits in the first variable in theorem \ref{chernbs} we obtain the following statement. 
\begin{theorem}\label{corhomch}
Let $ G $ be a totally disonnected group and let $ X $ be a proper $ G $-simplicial complex. Then there exists a 
natural transformation 
$$
K^G_*(X) \rightarrow H^G_*(X)
$$
which becomes an isomorphism after tensoring the left hand side with the complex 
numbers.
\end{theorem}
For totally disconnected groups the universal space for proper actions $ \underline{E}G $
can be chosen to be a $ G $-simplicial complex. Applying theorem \ref{corhomch} to 
this particular space yields the following result. 
\begin{theorem} Let $ G $ be a totally disonnected group. Then there exists a
Chern character 
$$
ch^G_*: K^G_*(\underline{E}G) \rightarrow H^G_*(\underline{E}G)
$$
from the topological $ K $-theory of $ G $ to the cosheaf homology of $ \underline{E}G $ 
which is an isomorphism after tensoring the left hand side with $ \mathbb{C} $. 
\end{theorem}

\section{The character of Baum and Schneider}\label{secbs}

In this section we compare the Chern character obtained in section \ref{secchern} 
with the bivariant equivariant Chern character for profinite groups defined by Baum and Schneider in \cite{BS}. 
The construction of Baum and Schneider is based on universal coefficient theorems for 
equivariant $ KK $-theory and bivariant equivariant cohomology.  
Roughly speaking, we show that in the case of profinite groups the equivariant Chern character $ ch^G_* $ reduces 
to the character defined by Baum and Schneider. 
In this way we obtain in particular a convenient description of $ ch^G_* $ for profinite groups. 
For simplicity we consider only finite $ G $-simplicial complexes in this section. \\
Let $ G $ be a profinite group. To every locally compact $ G $-space $ X $ one associates the Brylinski space 
\begin{equation*}
\hat{X} = \{(t,x) \in G \times X|\,t \cdot x = x\} \subset G \times X
\end{equation*}
which is again a $ G $-space by considering the action given by
$$
s \cdot (t,x) = (s t s^{-1}, s \cdot x)
$$ 
for $ s \in G $ and $ (t,x) \in \hat{X} $. 
Note that if $ G $ is a finite group one may view $ \hat{X} $ as the disjoint union of 
the fixed point sets $ X^t = \{x \in X|\, t\cdot x = x\} $ of elements $ t \in G $. \\
Now let $ X $ and $ Y $ be finite $ G $-simplicial complexes. The Chern character constructed by 
Baum and Schneider is of the form
\begin{align*}
\phi^G_*: KK^G_*(C(X), C(Y)) \rightarrow \SHom_G(K^*(\hat{X})_{\mathbb{C}}, K^*(\hat{Y})_{\mathbb{C}}) 
\rightarrow \SHom_G(H^*(\hat{X}), H^*(\hat{Y}))
\end{align*}
where the subscript $ \mathbb{C} $ stands for tensoring with $ \mathbb{C} $ over the integers and the 
second arrow is induced by the classical Chern character for $ K $-theory. 
Similarly, we may view $ ch^G_* $ as a map 
$$
ch^G_*: KK^G_*(C(X), C(Y)) \rightarrow HP^G_*(C^\infty(X), C^\infty(Y)) \cong 
\SHom_G(H^*(\hat{X}), H^*(\hat{Y}))
$$ 
using the universal coefficient formula for bivariant equivariant cohomology obtained in \cite{BS}. \\
In order to compare these transformations we consider first the case that $ X $ is a point. 
Let us review the construction of Baum and Schneider in this situation. Given a $ G $-equivariant vector 
bundle $ E $ over a compact $ G $-space $ Y $ one considers the pull-back bundle $ \pi^* E $ 
along the natural projection $ \pi: \hat{Y} \rightarrow Y $. The bundle $ \pi^* E $ is again 
$ G $-equivariant and carries in addition an action of the profinite completion $ \hat{\mathbb{Z}} $ 
of the integers. To describe this action note that the space of sections of $ \pi^* E $ is an $ \ayd $-module 
in a natural way. On the level of sections, the action of $ \hat{\mathbb{Z}} $ is determined by the natural 
automorphism $ T $ of this $ \ayd $-module. 
In fact pull-back along $ \pi $ yields a homomorphism
$$
\pi^*: K^0_G(Y)\rightarrow K^0_{\hat{\mathbb{Z}}}(\hat{Y})^G
$$
of equivariant $ K $-groups where the superscript $ G $ on the right hand side denotes taking the invariant part 
under the action of $ G $. Taking the trace of the element $ 1 \in \hat{\mathbb{Z}} $ 
in a representation of $ \hat{\mathbb{Z}} $ defines a ring homomorphism 
$$
\tau: R(\hat{\mathbb{Z}}) \rightarrow \mathbb{C}.
$$
Since the action of $ \hat{\mathbb{Z}} $ on $ \hat{Y} $ is trivial one obtains a homomorphism 
$$
\Pi:\xymatrix{K^0_G(Y)\; \ar@{->}[r]^{\pi^*} & K^0_{\hat{\mathbb{Z}}}(\hat{Y})^G_\mathbb{C} \ar@{->}[r]^{\!\!\!\!\!\!\!\!\!\!\!\!\!\!\!\!\!\!\!\cong} & 
(K^0(\hat{Y}) \otimes_{\mathbb{Z}} R(\hat{\mathbb{Z}})_\mathbb{C})^G \ar@{->}[r]^{\;\;\qquad \id \otimes \tau} & K^0(\hat{Y})_{\mathbb{C}}^G
 }
$$
where, as above, the subscript $ \mathbb{C} $ stands for tensoring with $ \mathbb{C} $. The transformation 
$ \phi_0 $ is obtained by composing the map $ \Pi $ with the map induced by the ordinary Chern character 
$ \ch: K^0(\hat{Y}) \rightarrow H^*(\hat{Y}) $. We shall assume that the latter is normalized as in \cite{BGV} in order 
to be compatible with cyclic homology. 
In the odd case the character of Baum and Schneider is defined using suspension. 
\begin{prop}\label{bscomp1}
Let $ G $ be a finite group and let $ Y $ be a finite $ G $-simplicial complex. The
equivariant Chern character 
$$ 
ch^G_*: K^G_*(C(Y)) \rightarrow \bigoplus_{j \in \mathbb{Z}} H^{* + 2j}(\hat{Y})^G
$$ 
agrees with the transformation $ \phi^G_* $ constructed by Baum and Schneider in this case. 
\end{prop} 
\proof The algebra $ C^\infty(Y) $ is a dense subalgebra of $ C(Y) $ which is closed under 
holomorphic functional calculus. As in the nonequivariant situation the natural inclusion $ C^\infty(Y) \rightarrow C(Y) $ induces an 
isomorphism $ K_0^G(C^\infty(Y)) \cong K^G_0(C(Y)) $. 
Consequently, to obtain the assertion in the even case it suffices to compare the characters $ ch^G_0(p) $ and $ \phi_0(p) $ 
of a $ G $-invariant idempotent $ p \in C^\infty(Y) \otimes \LH(V) $ 
where $ V $ is a finite dimensional Hilbert space equipped with a unitary representation $ \lambda $ of $ G $. 
Here the algebra $ \LH(V) $ of linear operators on $ V $ is equipped with the natural action $ s \cdot T = \lambda(s) T \lambda(s^{-1}) $ induced by $ \lambda $. 
Under the Chern character $ \ch^G_0 $ the element 
in $ HP^G_0(\mathbb{C}, C^\infty(Y)) $ corresponding to $ p $ is given by the cycle 
$$
\ch^G_0(p)(s) = \tr(\lambda(s) p) + \sum_{j = 1}^\infty (-1)^j \frac{(2j)!}{j!} \tr \biggl(\lambda(s) \biggl(p - \frac{1}{2}\biggr) (dpdp)^j\biggr)
$$
in the equivariant periodic cyclic complex of $ C^\infty(Y) $ where $ \tr $ denotes the ordinary trace. Here we use the analogue of 
theorem \ref{homotopyeq} for the periodic theory and work with noncommutative differential forms instead of the 
$ X $-complex in the second variable. The above formula for the Chern character appears also in the work of Block and Getzler 
on actions of compact Lie groups \cite{BG}. \\
Since $ G $ is finite the space $ \hat{Y} $ is the disjoint union of the fixed point sets $ Y^s $ for all elements $ s \in G $. 
The equivariant Hochschild-Kostant-Rosenberg theorem \cite{Voigtbs} shows that $ ch^G_0(p) $ is given by the class 
$$
\sum_{s \in G} \sum_{j = 0}^\infty \frac{(-1)^j}{j!} \tr(\lambda(s) p (dpdp)^j)_{|Y^s}
$$
in the even cohomology of $ \hat{Y} $. It is straightforward to check that this class is equal to the 
Chern character of the element $ \Pi(p) \in K^0(\hat{Y})_{\mathbb{C}} $. 
Note that the occurence of the operator $ \lambda(s) $ 
arises from the map $ \tau $ above. As a consequence, we see that the cohomology classes 
$ ch^G_0(p) $ and $ \phi_0(p) $ are equal. \\
In the order to treat the odd case consider the space $ S^1 $ with trivial $ G $-action and the split extension 
$$
\xymatrix{
C_0(Y \times (0,1))\; \ar@{>->}[r] & C(Y \times S^1) \ar@{->>}[r]^{\;\; \iota^*} & C(Y)
}
$$
of $ G $-$ C^* $-algebras where $ \iota^* $ is the transpose of the inclusion 
$ \iota: Y \rightarrow Y \times S^1 $ given by $ \iota(y) = (y,1) $. 
One obtains a corresponding commutative diagram 
$$
\xymatrix{
K^G_0(C_0(Y \times (0,1)))\; \ar@{>->}[r] \ar@{->}[d]^{\phi^G_0} & K^G_0(C(Y \times S^1)) \ar@{->>}[r] \ar@{->}[d]^{\phi^G_0} 
& K_0^G(C(Y)) \ar@{->}[d]^{\phi^G_0} \\
H^*_c(\hat{Y} \times (0,1))^G \; \ar@{>->}[r] & H^*(\hat{Y} \times S^1)^G \ar@{->>}[r] & H^*(\hat{Y})^G 
 }
$$
for the transformation defined by Baum and Schneider. In order to apply the character $ ch^G_0 $ in this situation we choose an 
equivariant triangulation of $ Y \times S^1 $ 
such that the inclusion $ \iota: Y \rightarrow Y \times S^1 $ becomes a map of $ G $-simplicial complexes. 
Similarly, we choose an equivariant triangulation of $ Y \times (0,1) $. Remark that the 
inclusion $ Y \times (0,1) \rightarrow Y \times S^1 $ will not be a simplicial map with respect to these triangulations. Still, 
it is easily verified that the analogue of the previous diagram for $ ch^G_0 $ 
is commutative. \\
According to our previous considerations the transformations $ \phi^G_0 $ and $ ch^G_0 $ 
agree for $ Y \times S^1 $ and $ Y $. It follows that $ \phi^G_0 $ and $ ch^G_0 $ agree also for the 
suspension $ Y \times (0,1) $ of $ Y $. As already mentioned above, the character $ \phi^G_1 $ is obtained 
from the even case using suspension. Assuming that $ \phi^G_1 $ is normalized
correctly it follows that $ ch^G_* $ and $ \phi^G_* $ agree in the odd case as well. \qed \\
One can verify directly that the assertion of proposition \ref{bscomp1} holds also for profinite groups. 
We shall obtain this as a special case of the following result. 
\begin{theorem}
Let $ G $ be a profinite group. On the category of finite $ G $-simplicial complexes the 
equivariant Chern character 
$$ 
ch^G_*: KK^G_*(C(X), C(Y)) \rightarrow \bigoplus_{j \in \mathbb{Z}} H^{* + 2j}_G(X,Y)
$$ 
agrees with the natural transformation $ \phi^G_* $ constructed by Baum and Schneider. 
\end{theorem} 
\proof Let us first assume that $ G $ is finite. As before we use the identification
$$
\bigoplus_{j \in \mathbb{Z}} H^{* + 2j}_G(X,Y) \cong \SHom_G(H^*(\hat{X}), H^*(\hat{Y}))
$$
for the bivariant equivariant cohomology of $ X $ and $ Y $. It suffices to show for every $ t \in G $ that the transformations 
$ ch^G_* $ and $ \phi^G_* $ agree after localization in $ t $ on the right hand side. 
If $ Z_t $ denotes the centralizer of $ t $ in $ G $ and $ f \in KK^G_*(C(X), C(Y)) $ is an arbitrary class we shall thus prove that 
the elements in $ \Hom_{Z_t}(H^*(X^t),  H^*(Y^t)) $ induced by $ ch^G_*(f) $ and 
$ \phi^G_*(f) $ are equal. \\
The localized Chern character $ ch^G_* $ can be written as the composition
\begin{align*}
ch^G_*: KK^G_*(C(X), C(Y)) \rightarrow KK^{Z_t}_*(C(X), C(Y))
\rightarrow \Hom_{Z_t}(H^*(X^t), H^*(Y^t))
\end{align*}
where the first arrow is the obvious forgetful map and the second arrow is the localized Chern character 
$ ch^{Z_t}_* $ for $ Z_t $. There is an analogous factorization 
for the localized character $ \phi^G_* $ due to its construction. 
Hence it suffices to consider the group $ G = Z_t $ and localization in $ t $. \\
Let us moreover specialize to $ X = Z_t/(t) $ where $ (t) $ denotes the subgroup of $ Z_t $ generated by $ t $. 
The cohomology $ H^*(Z_t/(t)) \cong C(Z_t/(t))[0] $ is located in degree zero 
and we have a commutative diagram 
$$
\xymatrix{
KK^{Z_t}_*(C(Z_t/(t)), C(Y)) \; \ar@{->}[r]^{\quad\; \cong} \ar@{->}[d]^{ch^{Z_t}_*} & KK^{(t)}_*(\mathbb{C}, C(Y)) \ar@{->}[d]^{ch^{(t)}_*} \\
\Hom_{Z_t}(C(Z_t/(t))[0], H^*(Y^t)) \ar@{->}[r]^{\quad\;\, \cong} & \Hom_{(t)}(\mathbb{C}[0], H^*(Y^t))
}
$$
obtained by induction. Note that $ \Hom_{(t)}(\mathbb{C}[0], H^*(Y^t)) = H^*(Y^t) $ since the action of $ t $ on $ Y^t $ is trivial. 
Inspecting the construction in \cite{BS} we see that there is a corresponding commutative diagram 
where the vertical maps are replaced by the transformations $ \phi^{Z_t}_* $ and $ \phi^{(t)}_* $, respectively. 
According to proposition \ref{bscomp1} the right vertical arrows in these diagrams coincide. It follows that the localized characters 
for $ Z_t $ are equal in the case $ X = Z_t/(t) $. \\
Both transformations $ \phi^{Z_t}_* $ and $ ch^{Z_t}_* $ are multiplicative with respect 
to the Kasparov product and the composition product. Let $ f \in KK^{Z_t}_*(C(X), C(Y)) $ 
be an arbitrary element and consider $ v \in KK^{Z_t}_*(C(Z_t/(t)), C(X)) $. 
After localization in $ t $ we obtain  
$$
ch^{Z_t}_*(v) \cdot \phi^{Z_t}_*(f) = \phi^{Z_t}_*(v) \cdot \phi^{Z_t}_*(f) = \phi^{Z_t}_*(v \cdot f) = 
ch^{Z_t}_*(v \cdot f) = ch^{Z_t}_*(v) \cdot ch^{Z_t}_*(f)
$$
according to our previous discussion. Since the images of elements of the form $ ch^{Z_t}_*(v) $ in 
$ \Hom_{(t)}(\mathbb{C}[0], H^*(X^t)) $ generate $ H^*(X^t) $ it follows that $ \phi^{Z_t}_*(f) $ and $ ch^{Z_t}_*(f) $ define the 
same element in $ \Hom_{Z_t}(H^*(X^t), H^*(Y^t)) $. We conclude that the characters $ \phi^G_* $ and $ ch^G_* $ agree 
in the case that $ G $ is finite. \\
For an arbitrary profinite group $ G $ the construction of $ \phi^G_* $ in section 3 of \cite{BS} shows that one 
has a commutative diagram 
$$
\xymatrix{
\varinjlim_{H \subset G} KK^{G/H}_*(C(X), C(Y/H))\; \ar@{->}[d]^{\varinjlim \phi^{G/H}_*} \ar@{->}[r]^{\qquad\quad \cong}  &  KK^G_*(C(X), C(Y)) \ar@{->}[d]^{\phi^G_*} \\
\varinjlim_{H \subset G} \bigoplus_{j \in \mathbb{Z}} H_{G/H}^{* + 2j}(X, Y/H) \;\ar@{->}[r]^{\qquad \cong} & \bigoplus_{j \in \mathbb{Z}} H_G^{* + 2j}(X, Y) 
}
$$
where the limits are taken over all open normal subgroups $ H $ of $ G $ which act trivially on $ X $. 
It is easily verified that there is an analogous commutative diagram where the vertical arrows are replaced by 
$ \varinjlim ch^{G/H}_* $ and $ ch^G_* $, respectively. According to our discussion for 
finite groups the left vertical maps in these diagrams are equal. Hence the same is true 
for the right vertical maps $ \phi^G_* $ and $ ch^G_* $. \qed

\section{The homological Chern character and equivariant Bredon homology} \label{seccw}

In this section we consider the homological Chern character otained by inserting the complex 
numbers in the second variable of the equivariant $ KK $-groups. The target of this 
transformation can be identified with equivariant Bredon homology. 
We show how to extend the homological Chern character
obtained in section \ref{secchern} appropriately from proper $ G $-simplicial complexes to 
proper $ G $-$ CW $-complexes. \\
Let $ G $ be a totally disconnected group. We call a $ G $-space $ X $ smooth if all isotropy groups $ G_x $ for 
$ x \in X $ are open subgroups of $ G $. A $ G $-space $ X $ is obtained from the 
$ G $-space $ A $ by attaching smooth equivariant $ k $-dimensional cells if there 
is a $ G $-equivariant pushout
$$
\xymatrix{
\coprod_{i \in I} G/H_i \times S^{k - 1} \; \ar@{->}[r] \ar@{->}[d] & A \ar@{->}[d] \\
\coprod_{i \in I} G/H_i \times D^k\;\ar@{->}[r] & X 
}
$$
where $ (H_i)_{i \in I} $ is a family of open subgroups of $ G $. 
A smooth $ G $-$ CW $-complex is a 
$ G $-space $ X $ together with a $ G $-invariant filtration 
$$
\emptyset = X_{-1} \subset X_0 \subset X_1 \subset \cdots \subset \bigcup_{k = 0}^\infty X_k = X 
$$
such that $ X $ carries the weak topology with respect to this filtration and for 
every $ k $ the space $ X_k $ is obtained from $ X_{k - 1} $ by attaching smooth equivariant $ k $-dimensional cells. 
In the sequel all $ G $-$ CW $-complexes are assumed to be smooth. This means in particular that every 
$ G $-$ CW $-complex can be viewed as an ordinary $ CW $-complex in a natural way. 
Note that every $ G $-simplicial complex is a $ G $-$ CW $-complex. 
A $ G $-$ CW $-complex $ X $ is called $ G $-finite if the quotient $ X/G $ is compact. 
It is proper iff all isotropy groups $ G_x $ for $ x \in X $ are compact open subgroups of $ G $.
For more detailed information about $ G $-$ CW $-complexes we refer to \cite{LueckLNM}. \\
We denote by $ sd^n(X) $ the $ n $th iterated barycentric subdivision of a $ G $-simplicial 
complex $ X $. Note that the barycentric subdivision of a $ G $-simplical complex is again a 
$ G $-simplicial complex. The following equivariant simplicial approximation theorem
is a straightforward extension of ordinary simplicial approximation as it can be found, for 
instance, in \cite{Spanier}. 
\begin{prop}\label{simpapprox} Let $ f: X \rightarrow Y $ be an equivariant continuous 
map between $ G $-simplicial complexes. If $ X $ is $ G $-finite there 
exists a natural number $ n $ and an equivariant simplicial map $ F: sd^n(X) \rightarrow Y $ 
equivariantly homotopic to $ f $. 
\end{prop} 
The next proposition contains the basic ingredient needed to extend the equivariant Chern character $ ch^G_* $ to 
$ G $-$ CW $-complexes. 
\begin{prop}\label{simpcw}
Let $ X $ be a $ G $-finite $ G $-$ CW $ complex. Then there exists a $ G $-finite 
$ G $-simplicial complex which is equivariantly homotopy equivalent to $ X $. 
\end{prop}
\proof We use induction on the dimension of $ X $, the case $ \dim(X) = 0 $ being clear. Using lemma 2.13 in \cite{Lueck} it 
suffices to show that the space obtained by attaching a finite number of smooth equivariant $ k $-dimensional cells 
to a $ (k - 1) $-dimensional $ G $-simplical complex $ Y $ is equivariantly homotopy equivalent to a 
$ G $-simplicial complex. Hence let us consider a $ G $-pushout diagram 
$$
\xymatrix{
\coprod_{i \in I} G/H_i \times \partial \Delta^k \; \ar@{->}[r]^{\qquad \quad f} \ar@{->}[d] & Y \ar@{->}[d] \\
\coprod_{i \in I} G/H_i \times \Delta^k \; \ar@{->}[r] & X 
}
$$
where $ \Delta^k $ denotes the standard $ k $-simplex and $ I $ is finite. 
According to proposition \ref{simpapprox} there exists $ n \in \mathbb{N} $ such that 
the upper horizontal map $ f $ is equivariantly homotopic to an equivariant simplicial map $ g: \coprod_{i \in I} G/H_i \times 
sd^n(\partial \Delta^k) \rightarrow Y $. 
Using again lemma 2.13 in \cite{Lueck} 
it follows that the $ G $-pushouts corresponding to $ f $ and $ g $ are equivariantly homotopy equivalent.
Let us abbreviate $ Y_0 = \coprod_{i \in I} G/H_i \times sd^n(\partial \Delta^k) $ 
and $ Y_1 = \coprod_{i \in I} G/H_i \times sd^n(\Delta^k) $ and consider the map $ g: Y_0 \rightarrow Y $. 
After replacing $ Y_0 $ by an equivariantly homotopy equivalent $ G $-simplicial complex $ Z_0 $ 
we can achieve that neighbouring vertices, that is, vertices which are connected by an edge, are mapped to different
vertices in $ Y $. If we also replace $ Y_1 $ appropriately by an equivariantly homotopy equivalent $ G $-simplicial 
complex $ Z_1 $ we obtain a $ G $-pushout diagram 
$$
\xymatrix{
Z_0 \; \ar@{->}[r]^h \ar@{->}[d] & Y \ar@{->}[d] \\
Z_1 \; \ar@{->}[r] & Z
}
$$
of $ G $-finite $ G $-simplicial complexes and equivariant simplicial maps. 
Moreover, by construction, in the commutative diagram 
$$
\xymatrix{
Y_1 \; \ar@{<-}[r] \ar@{->}[d] & Y_0 \ar@{->}[r]^g \ar@{->}[d] & Y \ar@{=}[d] \\
Z_1 \; \ar@{<-}[r] & Z_0 \ar@{->}[r]^h & Y  
}
$$
the left horizontal arrows are $ G $-cofibrations and the vertical maps are equivariant homotopy equivalences.
Using once more lemma 2.13 in \cite{Lueck} 
we conclude that $ Z $ and $ X $ are equivariantly homotopy equivalent. This proves the assertion. \qed \\ 
Let $ k^G $ and $ h^G $ be homotopy invariant functors on the 
category of $ G $-finite proper $ G $-$ CW $-complexes and equivariant maps. Moreover assume that 
$ \phi: k^G \rightarrow h^G $ is a natural transformation on 
the subcategory of $ G $-finite proper $ G $-simplicial complexes and equivariant 
simplicial maps. We shall show how $ \phi $ can be extended to a natural transformation on the whole 
category. This is similar to the techniques explained in the appendix of \cite{Dold}. \\
For a proper $ G $-finite $ G $-$ CW $-complex $ X $ let $ \mathcal{S}_X $ be the following 
category. The objects in $ \mathcal{S}_X $ are the equivariant homotopy equivalences 
$ S \rightarrow X $ where $ S $ is a proper $ G $-finite $ G $-simplicial complex. 
Note that there exist such equivariant homotopy equivalences due to proposition \ref{simpcw}. 
A morphism from $ R \rightarrow X $ to $ S \rightarrow X $ in $ \mathcal{S}_X $ is an equivariant 
simplicial map $ f: R \rightarrow S $ such that the diagram 
$$
\xymatrix{
R \; \ar@{->}[r] \ar@{->}[d]^f & X \ar@{=}[d] \\
S \; \ar@{->}[r] & X 
}
$$
is  commutative up to equivariant homotopy. To be precise, we shall rather work with 
isomorphism classes of $ G $-simplicial complexes in order to achieve that the category $ \mathcal{S}_X $ is small. 
Let us define 
$$
Sk^G(X) = \varinjlim_{S \rightarrow X} k^G(S)
$$
where the limit is taken over $ \mathcal{S}_X $. Using proposition \ref{simpapprox} 
one checks that $ Sk^G(X) $ is a functor on the category of $ G $-finite proper $ G $-$ CW $-complexes.
There is a canonical natural transformation $ Sk^G \rightarrow k^G $ which is an isomorphism 
since $ k^G $ was supposed to be homotopy invariant. In a similar way we obtain a functor 
$ Sh^G $ and a natural isomorphism $ Sh^G \rightarrow h^G $. 
For every $ S \rightarrow X $ in $ \mathcal{S}_X $ let us define 
$ \Phi_S(X): k^G(S) \rightarrow Sh^G(X) $ to be the composition
$ k^G(S) \rightarrow h^G(S) \rightarrow Sh^G(X) $ of $ \phi(S) $ with the canonical map. 
Passing to the limit we obtain a family of maps $ \Phi(X): Sk^G(X) \rightarrow Sh^G(X) $. 
Using again proposition \ref{simpapprox} this family is 
easily seen to be a natural transformation on the category of $ G $-finite proper $ G $-$ CW $-complexes. Due to 
the natural isomorphisms $ Sk^G \cong k^G $ and $ Sh^G \cong h^G $ 
we can view $ \Phi $ as a natural transformation from $ k^G $ to $ h^G $. By construction 
one obtains $ \Phi(X) = \phi(X) $ for every $ G $-simplicial complex $ X $. Moreover, the extended transformation $ \Phi $ is 
uniquely determined by $ \phi $. We have thus proved the following 
statement. 
\begin{prop}\label{homextension}
Let $ k^G $ and $ h^G $ be homotopy invariant functors on the category of $ G $-finite proper $ G $-$ CW $-complexes. 
Every natural transformation $ k^G \rightarrow h^G $ defined on 
the subcategory of $ G $-simplicial complexes and equivariant simplicial maps can 
be uniquely extended to a natural transformation $ k^G \rightarrow h^G $ on the whole category of $ G $-finite proper 
$ G $-$ CW $-complexes and equivariant continuous maps.
\end{prop}
Next we recall the definition of Bredon homology \cite{Lueck}, \cite{Voigtbs}. The smooth orbit category $ \Or(G) $ of a 
totally disconnected group $ G $ has as 
objects all homogenous spaces $ G/H $ where $ H $ is an 
open subgroup of $ G $. The morphisms in $ \Or(G) $ are all $ G $-equivariant maps. One obtains 
subcategories of $ \Or(G) $ by restricting the class of subgroups. 
We are interested in the class $ \mathcal{F} $ of all compact open subgroups of $ G $. 
The corresponding full subcategory $ \Or(G,\mathcal{F}) $ of $ \Or(G) $ consists of all 
homogeneous spaces $ G/H $ where $ H $ is compact open. \\
We will work over the complex numbers in the sequel. If $ \mathcal{C} $ is a small category a covariant 
(contravariant) $ \mathcal{C} $-vector space is a 
covariant (contravariant) functor from $ \mathcal{C} $ to the category of vector spaces. Morphisms 
of $ \mathcal{C} $-vector spaces are natural transformations. 
More generally one defines covariant and contravariant $ \mathcal{C} $-objects as 
functors with values in arbitrary target categories. \\
Given a contravariant $ \mathcal{C} $-vector space $ M $ and a covariant $ \mathcal{C} $-vector space 
$ N $ the tensor product $ M \otimes_\mathcal{C} N $ is the direct sum of $ M(c) \otimes N(c) $ over all 
objects $ c \in \mathcal{C} $ divided by the tensor relations $ mf \otimes n - m \otimes fn $ 
for $ m \in M(d), n \in N(c) $ and morphisms $ f: c \rightarrow d $ in $ \mathcal{C} $. \\
Let $ X $ be a proper $ G $-$ CW $-complex. There is
a contravariant functor from $ \Or(G, \mathcal{F}) $ to the category of $ CW $-complexes 
which associates to $ G/H $ the fixed point set $ X^H $. 
Composition with the covariant functor from $ CW $-complexes to chain complexes  
which associates to a $ CW $-complex $ Y $ the cellular chain complex $ C_*(Y) $ with complex coefficients
yields a contravariant $ \Or(G,\mathcal{F}) $-chain complex $ C_*^{\Or(G,\mathcal{F})}(X) $. \\
We define a covariant $ \Or(G, \mathcal{F}) $-vector space $ \mathcal{R}_q $ as follows. For a compact open subgroup $ H $ of $ G $ set 
$$ 
\mathcal{R}_q(G/H) = K_q(C^*(H)) \otimes_{\mathbb{Z}} \mathbb{C}
$$
where $ K_* $ denotes topological $ K $-theory and $ C^*(H) $ is the group $ C^* $-algebra of $ H $. Note that 
$ K_0(C^*(H)) = R(H) $ is the representation ring of $ H $ and $ K_1(C^*(H)) = 0 $. 
The character map induces an isomorphism 
$$ 
\mathcal{R}_0(G/H) = K_0(C^*(H)) \otimes_\mathbb{Z} \mathbb{C} \cong \mathcal{R}(H) 
$$ 
where $ \mathcal{R}(H) $ is the ring of conjugation invariant smooth functions on $ H $. \\
We define a chain complex $ C_*^{\Or(G,\mathcal{F})}(X;\mathcal{R}) $ by equipping 
$$
C_*^{\Or(G,\mathcal{F})}(X; \mathcal{R}) = \bigoplus_{p + q = *} C_p^{\Or(G,\mathcal{F})}(X) \otimes_{\Or(G,\mathcal{F})} \mathcal{R}_q
$$
with the differential induced from $ C_*^{\Or(G,\mathcal{F})}(X) $. 
\begin{definition}
Let $ G $ be a totally disconnected group. The equivariant Bredon homology of a proper $ G $-$ CW $-complex $ X $ (with 
coefficients in $ \mathcal{R}$) is 
$$
\mathcal{B}H^G_*(X; \mathcal{R}) = H_*(C_*^{\emph{\Or}(G,\mathcal{F})}(X; \mathcal{R})). 
$$ 
\end{definition}
The following result from \cite{Voigtbs} describes the relation between Bredon homology and 
cosheaf homology. Recall that in section \ref{secchern} the latter was defined in terms of equivariant cyclic homology. 
\begin{prop}\label{hpbredon}
Let $ G $ be a totally disconnected group and let $ X $ be a proper $ G $-simplicial complex. Then there 
is a natural isomorphism 
$$
H^G_*(X) \cong \mathcal{B}H^G_*(X; \mathcal{R})
$$ 
where $ H^G_* $ denotes cosheaf homology. 
\end{prop}
We shall now extend the homological Chern character obtained as a special case of theorem \ref{corhomch} to a natural transformation 
on the category of proper $ G $-CW-complexes. 
\begin{theorem}\label{cherncw}
Let $ X $ be a proper $ G $-$ CW $-complex. There exists a 
natural transformation 
$$
ch^G_*: K^G_*(X) \rightarrow \mathcal{B}H^G_*(X; \mathcal{R})
$$
which is an isomorphism after tensoring the left hand side with $ \mathbb{C} $. 
\end{theorem}
\proof Set $ k^G_* = K^G_* $ and let $ h^G_* $ be equivariant Bredon homology. 
According to proposition \ref{homextension} and proposition \ref{hpbredon}, the homological Chern character obtained 
in theorem \ref{corhomch} can 
be extended to the category of $ G $-finite proper $ G $-$ CW $-complexes. Since both $ k^G_* $ and 
$ h^G_* $ are theories with $ G $-compact supports, this transformation extends 
uniquely to the whole category of proper 
$ G $-$ CW $-complexes. Moreover theorem \ref{corhomch} implies that the extended transformation becomes an isomorphism
after tensoring the left hand side with the complex numbers. \qed \\
In \cite{Lueck} L\"uck has constructed Chern characters for proper equivariant homology theories in the context 
of discrete groups. A proper equivariant homology theory is an assignment which associates to 
every group $ G $ a proper $ G $-homology theory such that the theories for different groups are related 
by an induction structure. 
As a very special case of this general construction L\"uck obtains a natural isomorphism 
$$ 
\lambda_*^G: \mathcal{B}H^G_*(X; \mathcal{R}) \rightarrow K^G_*(X) \otimes_{\mathbb{Z}} \mathbb{C} 
$$ 
on the category of proper $ G $-$ CW $-complexes for a discrete group $ G $. 
\begin{prop}
Let $ G $ be a discrete group and let $ H $ be a finite subgroup. Then the map 
$$ 
\lambda_*^G: \mathcal{B}H^G_*(G/H; \mathcal{R}) \rightarrow K^G_*(G/H) \otimes_{\mathbb{Z}} \mathbb{C} 
$$ 
is inverse to $ ch^G_* $. 
\end{prop}
\proof The construction of L\"uck is compatible with the induction structures of the corresponding equivariant homology theories. 
In particular, induction from $ H $ to $ G $ yields a commutative diagram 
$$
\xymatrix{
\mathcal{B}H^G_*(G/H; \mathcal{R}) \; \ar@{->}[r]^{\lambda_*^G} \ar@{->}[d]^\cong & K^G_*(G/H) \otimes_{\mathbb{Z}} \mathbb{C} \ar@{->}[d]^\cong \\
\mathcal{B}H^H_*(\star\; ; \mathcal{R}) \; \ar@{->}[r]^{\lambda_*^H} & K^H_*(\star) \otimes_{\mathbb{Z}} \mathbb{C}
}
$$
where $ \star $ denotes the trivial $ H $-space consisting of a single point. Our Chern character 
is compatible with induction from $ H $ to $ G $ in the same way. 
Hence it suffices to consider a finite group $ H $ acting on a point. 
In this case the map $ \lambda_*^H $ is given by the identity map 
$ \mathcal{R}(H) \rightarrow \mathcal{R}(H) \cong R(H) \otimes_\mathbb{Z} \mathbb{C} $. 
On the other hand, the equivariant Chern character $ ch^H_* $ is induced by the character map $ R(H) \rightarrow \mathcal{R}(H) $ 
in this situation. This is precisely the map used for the identification $ R(H) \otimes_\mathbb{Z} \mathbb{C} \cong \mathcal{R}(H) $ which 
yields the assertion. \qed \\
It would be interesting to check how the transformations $ \lambda_*^G $ and $ ch^G_* $ between $ \mathcal{B}H^G_*(X; \mathcal{R}) $ and 
$ K^G_*(X) \otimes_{\mathbb{Z}} \mathbb{C} $ are related when $ X $ is an arbitrary proper $ G $-$ CW $-complex. 

\bibliographystyle{plain}

\end{document}